# HIERARCHICAL TESTING DESIGNS FOR PATTERN RECOGNITION

By Gilles Blanchard[1] and Donald Geman[2]

*CNRS and Fraunhofer FIRST, and Johns Hopkins University*

We explore the theoretical foundations of a "twenty questions" approach to pattern recognition. The object of the analysis is the computational process itself rather than probability distributions (Bayesian inference) or decision boundaries (statistical learning). Our formulation is motivated by applications to scene interpretation in which there are a great many possible explanations for the data, one ("background") is statistically dominant, and it is imperative to restrict intensive computation to genuinely ambiguous regions.

The focus here is then on pattern filtering: Given a large set $\mathcal{Y}$ of possible patterns or explanations, narrow down the true one $Y$ to a small (random) subset $\widehat{Y} \subset \mathcal{Y}$ of "detected" patterns to be subjected to further, more intense, processing. To this end, we consider a family of hypothesis tests for $Y \in A$ versus the nonspecific alternatives $Y \in A^c$. Each test has null type I error and the candidate sets $A \subset \mathcal{Y}$ are arranged in a hierarchy of nested partitions. These tests are then characterized by scope ($|A|$), power (or type II error) and algorithmic cost.

We consider sequential testing strategies in which decisions are made iteratively, based on past outcomes, about which test to perform next and when to stop testing. The set $\widehat{Y}$ is then taken to be the set of patterns that have not been ruled out by the tests performed. The total cost of a strategy is the sum of the "testing cost" and the "postprocessing cost" (proportional to $|\widehat{Y}|$) and the corresponding optimization problem is analyzed. As might be expected, under mild assumptions good designs for sequential testing strategies exhibit a steady progression from broad scope coupled with low power to high

Received March 2003; revised June 2004.

[1]Supported in part by a grant from the Humboldt Foundation and the IST Programme of the European Community under the PASCAL Network of Excellence, IST-2002-506778.

[2]Supported in part by ONR Contract N000120210053, ARO Grant DAAD19-02-1-0337 and NSF ITR DMS-02-19016.

power coupled with dedication to specific explanations. In the assumptions ensuring this property a key role is played by the ratio cost/power. These ideas are illustrated in the context of detecting rectangles amidst clutter.

**1. Introduction.** Motivated by problems in machine perception, specifically scene interpretation, we investigate the theoretical foundations of an approach to pattern recognition based on adaptive sequential testing. The basic scenario is familiar to everybody—identify one "pattern" (or "explanation") from among many by posing a sequence of subset questions. In other words, play a game of "twenty questions." Intuitively, we should ask more and more precise questions, progressing from general ones which "cover" many explanations, but are therefore not very selective, to those which are highly dedicated and decisive. Although the efficiency of coarse-to-fine (CTF) search drives the design of codes and many numerical routines, there has been surprisingly little work of a theoretical nature outside information theory to understand why this strategy is advantageous. We explore this question within the framework of sequential hypothesis testing, putting the emphasis on the modeling and optimization of computational cost: *In what sense and under what assumptions are the strategies which minimize total computation CTF?*

Needless to say, in order to have a feasible formulation of the problem one must make specific assumptions about the structure of the available tests (or "questions"). In this paper, we will therefore consider a particular structure based on an a priori multiresolution representation for the individual patterns and a corresponding hierarchy of hypothesis tests. Other important assumptions concern the statistical distribution of the tests and how cost varies with scope and power.

Our formulation is influenced by applications to pattern recognition, although we believe it remains sensible for other complex search tasks and we would argue that computational efficiency and CTF search are linked in a fundamental way. In both natural and artificial systems, many tasks do not require immediate, complete explanations of the input data. Nonetheless, the usual approach to machine perception is static: Intermediate results, when they exist, generally do not provide clear and useful provisional explanations. In contrast, we consider a sequence of increasingly precise interpretations (subsets of patterns), noting that experiments in biological vision (e.g., studies on "pop-out") report evidence for graded interpretations, for example, very fast identification of visual categories [27], "visual selection" and "regions of interest" [11].

Our formulation is also influenced by what we perceive to be some fundamental limitations in purely learning-based methods in pattern recognition in spite of recent advances (e.g., multiple classifiers, boosting and theoretical bounds on generalization error). We do not believe that very complex problems in machine perception, such as full-scale scene interpretation, will yield directly to improved methods of statistical learning. Some organizational framework is needed to confront the sheer number of explanations and complexity of the data (see, e.g., the discussion in [18]). In our approach learning comes into play in actually constructing the individual hypothesis tests from training data; in other words, one learns the individual components of an overall design.

The hypothesis-testing framework is as follows. Consider many patterns (or pattern classes) $y \in \mathcal{Y}$ as well as a special, dominating class 0 which represents "background." There is one true state $Y \in \{0\} \cup \mathcal{Y}$. In the highlighted applications, $Y$ refers to a semantic explanation of image data, for instance, the names and poses (geometrical presentations) of members belonging to a repertoire of actual objects appearing in an image. Thus, for example, a "pattern" might be a particular *instance* of a shape, say a square at some specific scale and orientation. The explanation $Y = 0$ represents "no pattern of interest" and is exceedingly more likely a priori; class 0 is also exceedingly more varied. Ultimately, we want to determine $Y$ (*classification* or *identification*). Ideally, this task would be accomplished rapidly and without error.

However, in machine perception and many other domains, near-perfect classification is often very difficult, even with sizable computational resources, and virtually impossible without resorting to a "contextual analysis" of competing explanations. In other words, we eventually need to test precise hypotheses $Y \in A$ against precise alternatives $Y \in B$, where $A, B \subset \mathcal{Y}$ ("*Is it an apple or a pear?*"). In view of the large number of possible explanations, it is not computationally feasible to anticipate all such scenarios. This argues for starting, and going as far as one can, with a "noncontextual analysis," meaning testing the hypothesis $Y \in A$ against the nonspecific alternative $Y \notin A$ (or, what is often almost the same, against the background alternative $Y = 0$) for a distinguished family of subsets $A \subset \mathcal{Y}$. Of course this only makes sense if there are *natural groupings* of explanations, which is certainly the case for pattern recognition (e.g., involving real objects and their spatial presentations).

Let $X_A$ denote the result of such a test, with $X_A = 1$ (resp. $X_A = 0$) indicating acceptance (resp. rejection). Indeed, it then makes sense to construct a family $\mathcal{X}$ of such tests *in advance*, say of order $\mathcal{O}(|\mathcal{Y}|)$. Throughout the paper we assume that the family $\mathcal{A}$ of sets $A \subset \mathcal{Y}$ for which (noncontextual) tests are built has a hierarchical, nested cell structure. These sets will be called *attributes* and their cardinality called their *scope*. In this scheme, the contextual analysis— testing against specific alternatives—begins only after the number of candidate explanations is greatly reduced, at which point tests may be created *on-line* to address the specific ambiguities encountered.

To pin things down, consider a toy example: Suppose $\mathcal{Y} = \{a, p\}$, standing for *apple* and *pear*, and $Y = 0$ stands for *other*, the most likely explanation. Suppose also there are four "tests":

(i) $X_{\{a,p\}}$ for testing $Y \in \{a, p\}$ versus $Y = 0$ (something like "*Is it a fruit ?*");

(ii) $X_{\{a\}}$ (resp. $X_{\{p\}}$) for testing $Y = a$ versus $Y = 0$ (resp. $Y = p$ vs. $Y = 0$);

(iii) $X_{\{avp\}}$ for testing $Y = a$ versus $Y = p$.

Tests $X_{\{a,p\}}, X_{\{a\}}, X_{\{p\}}$ are "noncontextual"; $X_{\{avp\}}$ is "contextual." Suppose all noncontextual tests have null false negative error. The type of CTF strategy that typically emerges from minimizing the "cost" of determining $Y$ under natural assumptions about how cost, scope and error are balanced is the intuitively obvious one: Perform $X_{\{a,p\}}$ first; then, if the result is positive ($X_{\{a,p\}} = 1$), perform $X_{\{a\}}$ and $X_{\{p\}}$; finally, perform $X_{\{avp\}}$ if both the previous results are again positive.

In this paper we consider efficient designs for the noncontextual phase only; the full problem, including contextual disambiguation, will be analyzed elsewhere. However, we anticipate the complexity of this contextual analysis by incorporating into our measurement of computation a "postprocessing" penalty which is proportional to the number of remaining explanations.

Our objective, then, is efficient "pattern filtering." The reduced set of explanations after noncontextual testing, denoted by $\widehat{Y}$ and called the set of *filtered patterns* (or *detected patterns*), is a *random subset* of $\mathcal{Y}$ that also depends on the chosen *strategy*, that is, the sequence of tests chosen to be performed. The tests are performed sequentially, and the choice of the next test to perform (or the decision to stop the search) depends on the outcomes of the past tests and is prescribed by the strategy. If strategy $T$ has performed the tests $X_{A_1}, \dots, X_{A_k}$ before terminating (note that $k$ and $A_2, \dots, A_k$ are themselves random variables), then the set of filtered patterns is determined in a simple way from the outcomes of the tests: $\widehat{Y}(T)$ consists of all patterns $y \in \mathcal{Y}$ which are "accepted" by every test $X_{A_i}$ for which $y \in A_i$, $1 \leq i \leq k$. In other words, a pattern is said to be filtered if it is not ruled out by one of the tests performed.

The fundamental constraint is no missed detections:

$$P(Y \in \widehat{Y} \cup \{0\}) = 1.$$

*This condition is satisfied if each individual test $X_A$ has zero type* I *error, and we make this assumption about every test $X_A$, recognizing that we must pay for it in terms of cost and power (or equivalently type* II *error). Although* we shall not be explicitly concerned with standard estimators such as

$$\widehat{Y}_{\mathrm{MLE}}(\mathcal{X}) = \arg\max_y P(\mathcal{X}|Y = y) \quad \text{and} \quad \widehat{Y}_{\mathrm{MAP}}(\mathcal{X}) = \arg\max_y P(Y = y|\mathcal{X}),$$

or even formulate a prior distribution for $Y$, it then follows that

$$P(\widehat{Y}_{\text{MLE}} \in \widehat{Y} \cup \{0\}) = P(\widehat{Y}_{\text{MAP}} \in \widehat{Y} \cup \{0\}) = 1.$$

Tests $X_A \in \mathcal{X}$ are then characterized by their scope ($|A|$), power $[P(X_A = 0 | Y \notin A)]$ and computational cost, and certain fundamental trade-offs are assumed to hold among these quantities. In order to accommodate differing applications and establish general principles, we will consider several scenarios, including both "fixed" and "variable" powers and two models—"power-based" and "usage-based"—for how the cost of a test is determined. Only the power-based cost model will be considered in detail; an analysis of the usage-based model can be found in [6]. Two other basic assumptions we make are (i) mean computation is well approximated by conditioning on $Y = 0$; and (ii), in that case, the tests are conditionally independent.

Except for a concluding illustration, we do not consider how these hypothesis tests $X_A$ are actually constructed, that is, depend functionally on the raw data. In the applications cited in Section 8 this typically involves statistical learning, for instance, inducing a decision tree or support vector machine from positive ($Y \in A$) and negative ($Y \notin A$) examples. *We are designing the specifications rather than the tests themselves, and modeling the computational process rather than learning decision boundaries for classification.* Presumably standard techniques can be used to build tests to the desired specifications if the trade-offs are reasonable. In Section 8 we will mention one recipe in an image analysis framework.

Although we will assume throughout that the true $Y$ is a single pattern belonging to $\{0\} \cup \mathcal{Y}$, our analysis would remain valid if we allowed $Y$ to be an entire subset of patterns $Y \subset \mathcal{Y}$ (with $Y = \varnothing$ representing "no pattern of interest" or "background"). In this case, $X_A$ would test the hypothesis $Y \cap A \neq \varnothing$ against $Y \cap A = \varnothing$, or against the nonspecific alternative $Y = \varnothing$. This setting might be more useful in some applications, such as scene interpretation, although in the end these subsets are simply more complex individual explanations.

Finally, our work is a natural outgrowth of an ongoing project on scene analysis (especially object recognition) which has been largely of an algorithmic nature (see, e.g., [2]). The current objective is to explore a suitable mathematical foundation. This was begun in [13] and [14] where the computational complexity of traversing abstract hierarchies was analyzed in the context of *purely* power-based cost—assuming that cost is an increasing, convex function of power. It was continued in [20], in which the optimality of depth-first CTF search for background-pattern separation [checking if $\widehat{Y}(\mathcal{X}) = \varnothing$] was established under the same model. The cost model here is more realistic because cost depends on scope as well as power.

Index of Main Notation

Objects:

$\mathcal{Y}$                set of all possible patterns or explanations

$Y \in \mathcal{Y} \cup \{0\}$        true (data-dependent) pattern (0 means background)

$P_0(\cdot)$                $= P(\cdot | Y = 0)$, the background distribution

Attributes:

$A$                    a grouping of objects (a.k.a. attribute)

$\mathcal{A}$                hierarchy of attributes

$\overline{\mathcal{A}}$                "augmented" hierarchy of attributes (see Section 4.5.4)

$\mathcal{Z}(\mathcal{A})$              coverings of $\mathcal{A}$: $\bigcup_{A \in Z} A = \mathcal{Y}$ for all $Z \in \mathcal{Z}(\mathcal{A})$

$A_1$                  coarsest attribute(s); root in the tree-structured case

Tests:

$X$                    binary random variable

$\beta(X) \in [0,1]$         $= P_0(X = 0)$, power of $X$

$c(X) \in [0, \infty)$        cost of $X$

$X_{A,\beta}$              test for attribute $A$ with power $\beta$

$\mathcal{X}$                family of tests indexed by $\mathcal{A}$; "fixed (powers) hierarchy"

$\widetilde{\mathcal{X}}$                family of tests indexed by $\mathcal{A}, \beta$; "variable-power hierarchy"

$\beta(A); c(A)$           power and cost of $X_A$ (fixed hierarchy case)

$\Gamma$                    increasing, subadditive complexity function for power-based cost

$\Psi$                    increasing, convex power function for power-based cost; $\Psi(0) = 0, \Psi(1) = 1$

Strategies:

$T$                    labeled binary tree, $T^\circ$ denotes internal nodes of $T$

$X(s); A(s); \beta(s)$       test at interior node $s$ of $T$; attribute and power of this test

$\mathcal{X}(t)$               set of tests along the branch leading to node $t$ of $T$

$\widehat{Y}(t) \subset \mathcal{Y}$         surviving (filtered) explanations at terminal node $t$ of $T$

$\widehat{Y}(T)$              filtered set of objects (surviving explanations) after testing

$q_X(T)$                probability of performing $X$ in $T$ under $P_0$

$C(T)$                 $= C_{\text{test}}(T) + C_{\text{post}}(\widehat{Y}(T))$: total cost

$C_{\text{test}}(T)$           random variable, sum of the costs of the tests performed in $T$

$C_{\text{post}}(\widehat{Y}(T))$         $= c^* |\widehat{Y}(T)|$, random variable, postprocessing cost

**2. Organization of the paper.** In Section 3 we provide a nontechnical overview of the results obtained in the paper. The precise mathematical setup appears in Section 4.

Our principal results appear in Sections 5–7. In Section 5, we consider the simplest case: There is one single test $X_A$ of fixed power and cost for each attribute $A \in \mathcal{A}$, and we present a fairly general sufficient condition under which CTF strategies are optimal. The "variable-power hierarchy" is examined in Section 6, namely a whole family of tests $(X_{A,\beta})$ for each attribute $A \in \mathcal{A}$ indexed by their power $\beta$. As the results for variable powers are decidedly not comprehensive, we attempt to strengthen the case for the "optimality" of CTF search with a variety of simulations at the end of Section 6. In Section 7, we mention a few analytical results for a substantially different cost model in which the cost of a test depends on the frequency with which it is used; this section amounts to a summary of results in [6].

In order to see how all this plays out in practice, we illustrate a few previous applications of this methodology to scene interpretation in Section 8. We also sketch an algorithm in Section 8 for a synthetic example of detecting rectangles in images against a background of "clutter"; the purpose is to illustrate in a controlled setting the quantities which figure in our analysis, especially how computation is measured and tests are constructed from data. Finally, in Section 9, we discuss some connections with related work and decision trees, critique our results and indicate some directions for future research.

**3. Overview of results.** A strategy $T$ can be represented as a binary tree with a test $X \in \mathcal{X}$ at each internal node and a subset $\widehat{Y}(t)$ at each external node or leaf $t$. The computational cost due to testing, $C_{\text{test}}(T)$, is a random variable—the sum of the costs of the tests performed before $\widehat{Y}$ is determined. The mean cost is then the average over all tests $X \in \mathcal{X}$ of the cost of $X$ weighted by the probability that $X$ is performed in $T$; these quantities will be defined more carefully in Section 4.

In anticipation of resolving the ambiguities in $\widehat{Y}$ in order to determine $Y$, we add to the mean testing cost a quantity which reflects the postprocessing cost, taken simply as $C_{\text{post}}(\widehat{Y}(T)) = c^*|\widehat{Y}(T)|$, where $c^*$ is a constant called the unit postprocessing cost. This charge may also be (formally) interpreted as the cost of performing perfect, albeit costly, tests for each individual nonbackground explanation in $\widehat{Y}$ in order to remove any remaining error under the background hypothesis [i.e., render $P(\widehat{Y} = \varnothing | Y = 0) = 1$]. The constant $c^*$ then represents the cost of a perfect individual test. Again, all tests have null false negative error, so "perfect" refers to full power.

The natural optimization question is then to find the strategy $T^*$ which minimizes the mean total computation:

$$T^* = \arg\min_T E[C(T)], \qquad C(T) = C_{\text{test}}(T) + C_{\text{post}}(\widehat{Y}(T)).$$

We are particularly interested in determining when $T^*$ is CTF in scope (meaning scope is decreasing along any root-to-leaf branch) and CTF in power (meaning power increases as scope decreases). Informally, the assumptions we impose are:

(a) *A multiresolution, nested cell representation*: The family of attributes $\mathcal{A}$ has the structure of a tree (see, e.g., Figure 1).

(b) *Background domination*: Mean computation $E[C(T)]$ and power $P(X_A = 0 | Y \notin A)$ are well approximated by taking $P = P_0 = P(\cdot | Y = 0)$.

(c) *Conditional independence*: Under $P_0$ families of tests over distinct attributes are independent. This is the strongest assumption and the one most likely to be violated in practice.

In the case of a fixed-powers hierarchy considered in Section 5, we assume that the test for attribute $A$ has cost $c(A)$ and power $\beta(A)$. We show that the ratios $c(A)/\beta(A)$ play a crucial role in the analysis of the optimization problem, and give the following general sufficient condition: *CTF optimality holds whenever, for any attribute $A$, the ratio of cost to power is less than the sum of the corresponding ratios over all direct children of $A$ in the test hierarchy* (including if necessary the perfect tests representing the postprocessing cost, having cost $c^*$ and power 1).

In the case of a variable-power hierarchy (Section 6), we consider a multiplicative model for the cost of $X_{A,\beta}$: $c(X_{A,\beta}) = \Gamma(|A|) \times \Psi(\beta)$, where $\Gamma$ is subadditive and $\Psi$ is convex. We prove that the CTF strategies always perform a specific test with the same power and that this power does not depend on the particular CTF strategy. A rigorous result about CTF optimality is only obtained for one particular $\Psi$, but simulations strongly indicate that the observed behavior is more widely true. In summary, CTF strategies seem to be optimal for a wide range of situations. The same can be said under the "usage-based" cost model in Section 7.

**4. Problem formulation.** In this section we formulate efficient pattern filtering as an appropriate optimization problem. (Recall that we are using the word "pattern" for an "explanation," often quite specific, rather than in the sense of some equivalence class of concepts or shapes.) We define the fundamental quantities which appear in this formulation, including attributes, tests and strategies, and how cost is measured both for individual tests and for testing designs. We also state our main assumptions about the test statistics and the relationships among cost, power and invariance which drive the optimization results in Sections 5–7.

4.1. *Goals.* The background probability space $\Omega$ represents the raw data—collections of numerical measurements—and $\mathcal{Y}$ denotes a set of *patterns* (or *classes* or *explanations*). We imagine the patterns $y \in \mathcal{Y}$ to be rather precise interpretations of the data and consequently $|\mathcal{Y}|$ to be very large. There is also a special explanation called *background*, denoted by 0, which represents "no pattern of interest" and is typically the most prevalent explanation by far.

We suppose there is a true state $Y$ which takes values in $\{0\} \cup \mathcal{Y}$ and which, for simplicity, is determined by the raw data. In other words, we regard $Y$ as a random variable on $\Omega$. Most of what follows could be generalized to the case in which $Y \subset \mathcal{Y}$ and $Y = \varnothing$ represents background.

EXAMPLE. In the context of machine perception, the raw data represent signals or images and the explanations represent the presentations of special entities, such as words in acoustical signals or physical objects in images (e.g., face instantiations or printed characters at a particular font and pose). The level of specificity of the explanations is problem-specific. However, we do assume that the data have in fact a unique interpretation at the level of precision of $Y$. Clearly this assumption eventually breaks down in the case of highly detailed semantic descriptions—at some point the subjectivity of the observer cannot be ignored.

The ultimate goal is *pattern identification*: Determine $Y$. However, for the reasons stated earlier, we shall focus instead on:

*Pattern filtering.* *Reduce the set of possible explanations to a relatively small, data-driven subset $\widehat{Y} \subset \mathcal{Y}$ such that $Y \in \widehat{Y} \cup \{0\}$ with probability (almost) 1.*

We shall also consider the special case of spotting one single, fixed pattern $y^*$. A related problem of interest is *background-pattern separation* or *background filtering*: Determine whether or not $Y = 0$. Background-pattern separation will not be analyzed in this paper since it has been studied elsewhere in a very similar framework in [13, 14, 20]. In contrast, detecting a single pattern of interest will often serve as a first step before turning to the filtering of all possible patterns. Formally, what will distinguish these tasks is only the postprocessing cost; see Section 4.5.2.

As discussed earlier, the rationale behind pattern filtering is that requiring that $Y \in \widehat{Y} \cup \{0\}$ ensures, by definition, that no pattern is missed. Hence, the ensuing analysis, which is aimed at determining $Y$ with high precision and is likely to be computationally intensive, can be limited to $\widehat{Y}$. Additional computation might involve a *contextual analysis*, such as constructing hypothesis tests on the fly for distinguishing between competing alternatives belonging to $\widehat{Y}$. This "postprocessing stage" will not be analyzed in this paper, *except* that we shall explicitly anticipate additional computation in the form of a penalty for unfinished business: We impose a "postprocessing cost" $C_{\text{post}}(\widehat{Y})$ proportional to the size of $\widehat{Y}$. The goal then is to find an optimal trade-off between the costs of "testing" and "postprocessing."

4.2. *Attributes and attribute tests.*    Any subset of patterns $A \subset \mathcal{Y}$ can be regarded as an "interpretation" of the data and we assume there are certain "natural groupings" of this nature (e.g., "writer" in a "Guess Who" version of twenty questions, "noun" in speech recognition and "character" in visual recognition). We call these distinguished subsets *attributes* and we denote the family of attributes by $\mathcal{A}$ (a collection of subsets of $\mathcal{Y}$) and suppose $|\mathcal{A}|$ is of order $\mathcal{O}(|\mathcal{Y}|)$. For every $y \in \mathcal{Y}$, we will assume that

$$(1) \qquad\qquad \{y\} = \bigcap_{A \ni y} A.$$

One of our main assumptions is that $\mathcal{A}$ has a multiresolution, hierarchical structure with attributes at varying levels of precision. Formally, we assume that

$$\forall A, A' \in \mathcal{A}, \qquad A \cap A' \neq \varnothing \Rightarrow (A' \subset A) \quad \text{or} \quad (A \subset A').$$

Note that the set of attributes thus has a tree structure (see Figure 1 for an example). Furthermore, assumption (1) implies that the set of leaves of the corresponding tree is exactly the set of all singleton attributes.

For every attribute $A \in \mathcal{A}$ we can build one or more binary *tests* $X$—the result of testing the hypothesis $Y \in A$ against either $Y \notin A$ or $Y = 0$; the value $X = 1$ corresponds to choosing $Y \in A$ and $X = 0$ to choosing the alternative. Which alternative, $Y \notin A$ or $Y = 0$, is more appropriate is application-dependent. For example, in inductive learning, the two cases correspond to the nature of the "negative" examples in the training set—whether they represent a random sample under $Y \in A^c$ or under "background." In the applications cited in Section 8, the tests are constructed based *entirely* on the statistical properties of the patterns in $A$; neither alternative is explicitly represented. Due to the domination of the background class, at least at the beginning of the search, and due to the simplification afforded by measuring total computation cost under $P_0 = P(\cdot | Y = 0)$, the

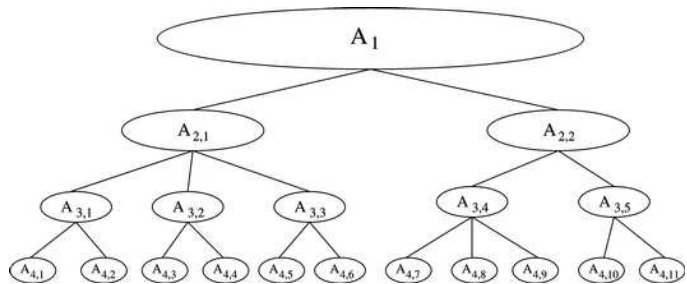

Fig. 1.    *Example of a (nonregular) tree-structured hierarchy of attributes.*

alternative hypothesis will hereafter be $Y = 0$ and we define the power of the test accordingly:

$$\beta(X) = P(X = 0 | Y = 0).$$

In order to make the notation more informative, we shall write either $X_A$ to indicate the attribute being tested or $X_{A,\beta}$ to signal both the attribute and the power.

The first main assumption we will make about these tests is that *their false negative rate is negligible.* In other words, if a pattern (i.e., nonbackground explanation) is present, then any attribute test which covers this pattern must respond positively:

(2) $$P(X_A = 1 | Y \in A) = 1 \qquad \forall A \in \mathcal{A}.$$

For this reason, and due to the origins of this work in visual object recognition, we sometimes refer to the size of $A$ as the *level of invariance* of $X_{A,\beta}$, but usually just as the *scope*. Its depth in the attribute hierarchy is called *level of resolution*. In general, however,

$$P(X_A = 1 | Y = y) > 0 \qquad \forall y \in \{0\} \cup \mathcal{Y} \setminus A.$$

In other words, the tests are usually not *perfect* or *two-sided invariants*.

Formally, assumption (2) is not necessary for the mathematical results in the coming sections to hold, because we will only make computations under the "background probability" (when $Y = 0$); see Section 4.4. However, this assumption is necessary for our formulation of pattern identification to make sense; indeed, it implies that if one has performed tests $X_{A_1}, \ldots, X_{A_k}$, then necessarily

$$Y \in \mathcal{Y} \setminus \bigcup_{k \,:\, X_{A_k} = 0} A_k.$$

We will say that the patterns above have been *filtered* by tests $X_{A_1}, \ldots, X_{A_k}$ and focus on sequential testing designs for which the chosen $\widehat{Y}$ is the set of patterns filtered by all the tests actually performed (called *filtered patterns* for short). This choice coheres with our requirement that $Y \in \widehat{Y} \cup \{0\}$ with probability 1, while at the same time ensuring that $\widehat{Y}$ is of minimum size given the available information.

Finally, each test $X_{A,\beta}$ has a *cost* or *complexity* $c(X_{A,\beta})$ which represents the amount of online computation (or time) necessary to evaluate $X_{A,\beta}$. In Section 4.6 we shall consider a cost model in which the cost of a test is a predetermined quantity related to power and scope. In Section 7 we briefly consider another "usage-based" cost model.

4.3. *Test hierarchies.* We consider two types of families of tests, one with exactly one test (at some fixed power) per attribute and referred to as a *fixed test hierarchy*, and one with a one-parameter family of tests $\{X_{A,\beta}, 0 \le \beta \le 1\}$ for each $A \in \mathcal{A}$ indexed by power and referred to as a *variable-power hierarchy*.

4.3.1. *Fixed hierarchy.* We will denote such a hierarchy by $\mathcal{X} = \{X_A, A \in \mathcal{A}\}$ and write $\beta(A)$ for the power of $X_A$ and $c(A)$ for its cost. Optimal testing strategies for fixed hierarchies is the subject of Section 5 and Section 7 for two different cost models. In the analysis in those sections a central role is played by the (random) set $\widehat{Y}(\mathcal{X})$ of patterns which are filtered by *all* the tests in $\mathcal{X}$, that is, those patterns which are verified at all levels of resolution. More precisely:

$$\widehat{Y}(\mathcal{X}) = \mathcal{Y} \setminus \bigcup \{A \in \mathcal{A} | X_A = 0\}.$$

Recall that under our constraint on the false negative error, we necessarily have $P(Y \in \{0\} \cup \widehat{Y}(\mathcal{X})) = 1$. Clearly, $\widehat{Y}(\mathcal{X})$ leads to a smaller postprocessing cost than any $\widehat{Y}$ based on only *some* of the tests in the hierarchy, but, of course, requires more computation to evaluate in general.

4.3.2. *Variable-power hierarchy.* The variable-power hierarchy is

$$\widetilde{\mathcal{X}} = \{X_{A,\beta} | A \in \mathcal{A}, \beta \in [0,1]\}.$$

In Section 6 we will consider testing strategies in which, at each step in a sequential procedure, *both* an attribute *and* a power may be selected. This clearly leads to a more complex optimization problem and our results in this direction are correspondingly far less complete than those in the case of a fixed hierarchy. From another point of view, extracting a subset of tests from a variable-power hierarchy (e.g., specifying a testing strategy) is a type of *model selection* problem.

4.4. *The probabilistic model.* In order for the upcoming optimization problems to be well defined, we need to specify the joint distribution of the random variables in $\widetilde{\mathcal{X}}$.

The first hypothesis we make is that we are going to measure mean computation relative to $P_0(\cdot) = P(\cdot | Y = 0)$—the "background distribution." This is justified by the assumption that a priori the probability of the explanation $Y = 0$ is far greater than the compound alternative $Y \ne 0$ let alone any single, nonbackground explanation. For instance, in visual processing a randomly selected subimage is very unlikely to support a precise explanation in terms of visible patterns; in other words, most of the time all we observe is clutter.

The second hypothesis we make is that, under $P_0$, any family of tests $X_{A_1,\beta_1}, \ldots, X_{A_k,\beta_k}$ for *distinct* attributes $A_1, \ldots, A_k$ is independent. This is probably the strongest assumption in this paper but is not altogether unreasonable under $P_0$ in view of the structure of $\mathcal{A}$ since two distinct tests are either testing for disjoint attributes (if they are at the same level of resolution) or testing for attributes at different levels of resolution. In Section 5 we shall briefly consider simulations for a nontrivial dependency structure—a Markov hierarchy.

No assumptions are made about the dependency structure among tests for the same attribute but at different powers. Instead, the assumption to be made in the following section that no attribute can be tested twice in the same procedure allows us to compare the cost of testing strategies regardless of this dependency structure.

4.5. *Testing strategies and their cost.* We consider sequential testing processes, where tests are performed one after another and the choice of the next test to be performed (or the decision to stop the testing process) can depend on the outcomes of the previously performed tests. We will make the important assumption that in any sequence of tests, a given attribute can only be tested once.

DEFINITION 1 (*Testing strategy*). A strategy is a finite labeled binary tree $T$ where each internal node $t \in T^\circ$ is labeled by a test $X(t) = X_{A(t),\beta(t)}$ and where $A(t) \neq A(s)$ for any two nodes $t, s$ along the same branch. At each internal node $t$ the right branch corresponds to $X(t) = 1$ and the left branch to $X(t) = 0$.

The restriction to at most one test per attribute $A$ along any given branch, whereas of course automatically satisfied in the case of a fixed hierarchy (Sections 5 and 7), does limit the set of possible strategies for a variable-power hierarchy since several tests $X_{A,\beta}$ of varying power are available for each attribute $A$. In that case the purpose of this assumption is essentially to simplify the analysis by guaranteeing that all the tests actually performed are independent.

The leaves (terminal nodes) of $T$ will be labeled in accordance with the answers to the tests: Every leaf of $T$ is labeled by the subset $\widehat{Y} \subset \mathcal{Y}$ of filtered patterns that have not been ruled out by the tests performed by the strategy (along the branch leading to this leaf). In other words, for any strategy $T$ and leaf $s$ of $T$, if $\mathcal{X}(s)$ denotes the set of tests along the branch leading to $s$, we put

$$\widehat{Y}(s) = \mathcal{Y} \setminus \bigcup \{A \in \mathcal{A} | X_A \in \mathcal{X}(s); X_A = 0\}.$$

The random set $\widehat{Y}(T)$ is then defined by interpreting $T$ as a function of the tests which takes values among its leaves. However, how the leaves are labeled is irrelevant for the purposes of defining the testing cost $C_{\text{test}}(T)$ of a strategy; it will only influence the postprocessing cost $C_{\text{post}}(\widehat{Y})$.

4.5.1. *Cost of testing.* There are several equivalent definitions of the testing cost of $T$, another random variable. One is

$$C_{\text{test}}(T) = \sum_{t \in T^\circ} c(X(t)) \mathbf{1}_{H_t},$$

where $H_t$ is the history of node $t$—the event that $t$ is reached. Recall that $T^\circ$ is the set of internal nodes of $T$. This is clearly the same as aggregating the costs over the branch traversed or adding the costs of all tests performed. Given a probability distribution $P$ on $\Omega$, and in particular $P = P_0$, two equivalent expressions for the *mean cost* are then

$$(3) \qquad E_0[C_{\text{test}}(T)] = \sum_{t \in T^\circ} c(X(t)) P_0(H_t) = \sum_X c(X) q_X(T),$$

where

$$q_X(T) = P_0(X \text{ performed in } T) = \sum_{t \in T^\circ} \mathbf{1}_{\{X(t)=X\}} P_0(H_t).$$

Expression (3) is particularly useful in proving some of our results; in Section 5 we will transform it into yet another expression that will anchor the analysis there.

4.5.2. *Cost of postprocessing.* It is natural to define the postprocessing cost in the following, goal-dependent manner:

(i) *Filtering a special pattern*: $C_{\text{post}}(\widehat{Y}(T)) = c^* \mathbf{1}_{\{y^* \in \widehat{Y}(T)\}}$ where $y^*$ is the target pattern.

(ii) *Filtering all patterns*: $C_{\text{post}}(\widehat{Y}(T)) = c^* |\widehat{Y}(T)|$.

Here $c^*$ is some constant called the *unit postprocessing cost.*

In the case of a single target pattern, note that this choice of postprocessing cost naturally leads us to disregard any attribute not containing the target $y^*$ as those tests are irrelevant to the goal at hand and can only augment the total cost. Consequently, the set of relevant attributes reduces to a "vine" $A_1 \supset A_2 \supset \cdots \supset A_L$. In this case, choosing a testing strategy boils down to choosing a subset of these relevant attributes and an order in which to test for them. If a test returns a null answer, the search terminates with the outcome $y^* \notin \widehat{Y}$ and there is no postprocessing charge; on the other hand, if all the selected tests respond positively, then $y^* \in \widehat{Y}$ is declared (which still may not be true) and the charge is $c^*$. In particular, the testing strategy $T$ itself has in this case the structure of a vine (see Figure 2). In contrast, in the case of general pattern filtering the testing strategies are of course tree-structured.

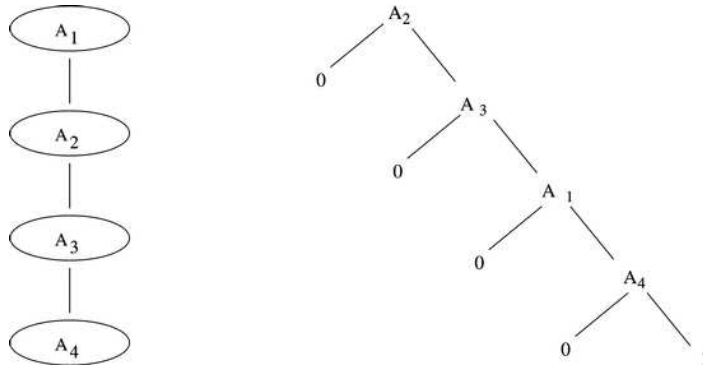

Fig. 2. Left: *A vine-structured hierarchy of attributes for detecting one pattern.* Right: *An example of a vine-structured testing strategy for this hierarchy.*

4.5.3. *Optimization problem.* The *total computational cost* for the task at hand is $C_{\text{test}}(T) + C_{\text{post}}(\widehat{Y}(T))$. The corresponding optimization problem, our central focus, is then to find a strategy attaining

$$(4) \qquad \min_{T \in \mathcal{T}} (E_0[C_{\text{test}}(T)] + E_0[C_{\text{post}}(\widehat{Y}(T))]),$$

where $\mathcal{T}$ is the family of all strategies. We emphasize that in the case of variable-power hierarchies we are therefore optimizing over both power and scope.

4.5.4. *Equivalent model with perfect tests.* There is an equivalent way to interpret the postprocessing cost which is technically more convenient. We can think of $c^*$ as the cost of performing a *perfect* test (i.e., without errors under $P_0$) for any individual pattern. Therefore, the postprocessing cost model is *formally equivalent* to supposing there is no postprocessing stage, but that *no errors* (under $P_0$) are allowed at the end of the procedure, enforced by performing, as needed, some additional perfect tests at the end of the search. Since we have assumed that no attribute, and in particular no singleton $\{y\}$, cannot be tested at two different powers along the same branch, we can incorporate perfect testing into the previous framework simply by adding a final layer to the original hierarchy $\mathcal{A}$ which copies the original leaves, thereby accommodating a battery of perfect singleton tests having cost $c^*$. (Conditional independence is actually maintained since the new tests are deterministic under $P_0$.) We denote by $\overline{\mathcal{A}}$ the resulting augmented hierarchy. [Due to this augmentation there is a slight abuse of notation when identifying an attribute with a subset of $\mathcal{Y}$, since in the augmented hierarchy we would like (in order to be entirely consistent) to consider some attributes as distinct although they correspond to the same set $\{y\}$. However, we will stick to the notation introduced before in order to avoid cumbersome changes.]

This formal construction allows us to include the postprocessing cost in the testing framework. Furthermore, in the augmented model it is not difficult to show that for any strategy $T$ there exists a strategy $T'$ performing exactly the same tests, but with the perfect tests performed at the end only, so that the optimization problem is in fact unchanged by allowing the perfect tests to be performed at any time. *In summary, the equivalent optimization problem is to minimize the amount of computation necessary to achieve no error under $P_0$ based on the augmented hierarchy.*

4.6. *Cost of a test.* There are certain natural trade-offs among cost, power and invariance:

(a) At a given cost, power should be a *decreasing* function of invariance.
(b) At a given power, cost should be an *increasing* function of invariance.
(c) At a given invariance, cost should be an *increasing* function of power.

In Section 5, we will first deal with a generic setting where the test associated to a given attribute $A$ has power $\beta(A)$ and cost $c(A)$. In Section 6 we will use a more specific model reflecting the trade-offs among cost, power and invariance mentioned above:

$$(5) \qquad\qquad c(X_{A,\beta}) = \Gamma(|A|) \times \Psi(\beta),$$

where the *complexity function* $\Gamma$ is subadditive and the *power function* $\Psi$ is convex. Consequently, we evaluate the cost of a test much like the merit of a dive in the Olympics: at any given level of difficulty ($\Gamma$) a score ($\Psi$) is assigned based on performance alone. For normalization, we can assume that $\Gamma(1) = 1$. Then with the equivalent model where the postprocessing cost is replaced by "perfect" tests in mind, it is consistent to assume $c^* = \Psi(1)$. *This multiplicative model is supported* (*at least roughly*) *by what is observed in actual experiments* (*see Section 8*).

One special case, treated in Section 6, is $\Gamma(n) = n$, that is, the complexity is simply the level of invariance. This case is the *least* favorable to CTF strategies since, in effect, no "credit" is given for shared properties among two disjoint attributes $A, B \in \mathcal{A}$. If, for instance, $|A| = |B|$ with $A, B$ disjoint, a test for $A \cup B$ at a given power $\beta$ has the same cost as testing separately for both $A$ or $B$ at power $\beta$.

A particular case, treated in [13] and [20], in the setting of a fixed hierarchy, is to assume $c(X_{A,\beta_A}) = \Psi(\beta_A)$ for some function $\Psi$. The model considered here is more general.

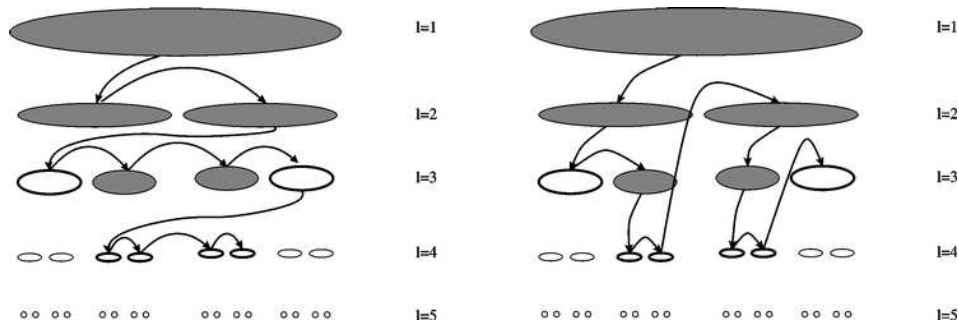

Fig. 3. *Example of typical CTF strategies.* Left: *breadth-first;* right: *depth-first.*

4.7. *Special strategies.* In the following sections our main goal will be to determine under what additional hypotheses the optimal strategies are "coarse-to-fine" (CTF).

DEFINITION 2 (*Coarse-to-fine*). A strategy $T \in \mathcal{T}$ is *CTF in resolution* or just *CTF*, if an attribute is tested if and only if each of its ancestors has already been tested and returned a positive answer. A strategy $T \in \mathcal{T}$ is *CTF in power* if, for any two nodes $s, t$ along the same branch, $\beta(s) \geq \beta(t)$ whenever $A(s) \subset A(t)$.

In the case of filtering a single pattern, this simply means that a CTF strategy performs all the relevant tests in the order of increasing resolution, that is, $X_{A_1}, \ldots, X_{A_L}$. For general pattern filtering, several different strategies have the CTF property, for instance, "breadth-first" and "depth-first" search. In Figure 3 these two CTF strategies are illustrated in the case of a hierarchy of depth $L = 5$ and test outcomes such that $\widehat{Y}(T_{\mathrm{ctf}}) = \varnothing$, that is, no patterns are verified at all resolutions due to the "null covering" $\{X_{3,1} = 0, X_{4,3} = 0, X_{4,4} = 0, \ X_{4,5} = 0, \ X_{4,6} = 0, X_{3,4} = 0\}$ (writing $X_{l,k}$ for the $k$th test at depth $l$). Notice that the *breadth-first* CTF strategy has the nice feature that the tests are always performed in the order of nondecreasing depth.

For a fixed hierarchy, all CTF strategies for pattern filtering perform exactly the same tests (although perhaps not in the same order). Whatever the order chosen, in the end, along any branch of the attribute hierarchy, every test has been performed starting from the root until the first null answer encountered on this branch. It is therefore possible to speak of "the" CTF strategy, it being understood that the precise order in which the tests are performed does not affect the mean cost.

NOTE. Whereas we do not consider the problem of separating patterns from background in and of itself (as in [20]), it is interesting to observe that

the situation is more complex in that case since all CTF strategies are not equivalent. Indeed, in any optimal strategy, testing stops as soon as any complete "1-chain" is found and, consequently, depth-first CTF strategies are generally optimal, as shown in [20].

The probability of performing a test $X_A$ in a CTF strategy for a fixed hierarchy has a simple expression:

$$q_A(T) = P_0(X_A \text{ performed in } T) = P_0(X_B = 1, \text{ for all } B \supset A, B \neq A)$$
$$= \prod_{B \supset A}(1 - \beta(B)).$$

Moreover, under the CTF strategy $\widehat{Y}$ minimizes $C_{\text{post}}(\widehat{Y}(T))$, and, in fact,

(6) $$\widehat{Y}(T_{\text{ctf}}) = \widehat{Y}(\mathcal{X}) \qquad \text{a.s.,}$$

which can also be identified with the set of all "1-chains" in the hierarchy. It follows that the total mean cost of the CTF strategy is then given by

$$E_0[C(T_{\text{ctf}})] = \sum_{A \in \mathcal{A}} c(A) \prod_{B \supset A}(1 - \beta(B)) + E_0[|\widehat{Y}(\mathcal{X})|].$$

Still in the case of a fixed hierarchy, it will be useful to delineate all strategies with property (6).

DEFINITION 3 (*Complete strategies*). A strategy $T \in \mathcal{T}$ is complete if $\widehat{Y}(T) = \widehat{Y}(\mathcal{X})$. The family of complete strategies is denoted by $\overline{\mathcal{T}}$.

REMARK. Under the hypotheses we have made, *for a complete strategy* it is possible to compute explicitly the probability of error under the null hypothesis before the postprocessing step, that is, to calculate $P_0(\widehat{Y} \neq \varnothing)$. (This is the probability that at least one nonnull pattern is detected when only background clutter is actually observed.) For single-pattern detection it is just the probability under $P_0$ that all the tests along the vine respond positively: $P_0(\widehat{Y} \neq \varnothing) = \prod_{k=1}^{L}(1 - \beta_k)$ [where $\beta_k = \beta(A_k) = P(X_{A_k} = 0)$]; for detection of all possible patterns it is exactly the probability that there exists a "1-chain" leading from the root of the attribute tree to one of its leaves. Given the independence assumption on the tests under $P_0$, this in turn is *exactly the probability of nonextinction of an inhomogeneous branching process at generation $L$, which can be computed explicitly once the branching probabilities* [*i.e.,* $\beta(A), A \in \mathcal{A}$] *are known.*

Finally, for a variable-power hierarchy $\widetilde{\mathcal{X}}$ there are many different, noncost-equivalent CTF strategies depending on the powers chosen for the tests along each branch. Nonetheless, surprisingly, the optimal CTF strategy can sometimes be precisely characterized, being CTF in power with, in fact, a unique power assigned to each attribute (see Section 6).

**5. Optimal strategies for fixed costs and powers.** Throughout this section we assume a fixed test hierarchy $\mathcal{X} = \{X_A, A \in \mathcal{A}\}$ and we write $c(A), \beta(A)$ for the cost and power, respectively, of $X_A$. We will then refer to "testing an attribute $A$" or "attaching an attribute" to a node of $T$ without ambiguity. Our goal is to identify conditions (trade-offs) involving $\{c(A), \beta(A), A \in \mathcal{A}\}$ under which optimal strategies may be characterized.

For parts of this section it will be easier to actually consider the equivalent model with perfect tests in lieu of the postprocessing cost, as described in Section 4.5.4. From here on, $\overline{\mathcal{A}}$ will denote the augmented hierarchy, and the considered strategies $T$ for $\overline{\mathcal{A}}$ will satisfy the no-error constraint. In other words, in the augmented model, when the strategy ends all patterns $y \in \mathcal{Y}$ must have been covered by at least one test which has been performed and returned 0 (again, it may be one of the perfect, artificial tests representing postprocessing). We start this section with a fundamental formula for the average cost $E_0[C(T)]$ that will be useful for all of the results to follow.

5.1. *Reformulation of the cost.* As just pointed out, in the augmented hierarchy model strategies must find a way to "cover" all patterns with attributes whose associated test is negative. Therefore the notion of *covering* will play a central role in the analysis to come, motivating the following definitions:

DEFINITION 4 (*Covering*). A set of attributes $\mathcal{Z} \subset \overline{\mathcal{A}}$ is a *covering* if

$$\bigcup \{A, A \in \mathcal{Z}\} = \mathcal{Y}.$$

The set of coverings for the augmented hierarchy $\overline{\mathcal{A}}$ is denoted $\mathcal{Z}(\overline{\mathcal{A}})$.

DEFINITION 5 (*Tested attributes*). For a given strategy $T$, denote by $\mathcal{X}(T)$ the (random) set of attributes tested by $T$, and by $\mathcal{X}_0(T)$ the set of attributes in $\mathcal{X}(T)$ for which the corresponding test returned the answer 0, called the *zero set* of $T$.

Of course, the no-error constraint for a strategy $T$ now reads simply: $\mathcal{X}_0(T)$ is (a.s.) a covering. We now turn to an important formula:

LEMMA 1 (*Cost reformulation*). *For any (no-error) strategy $T$ for the augmented hierarchy $\overline{\mathcal{A}}$,*

$$(7) \qquad E_0[C(T)] = \sum_{Z \in \mathcal{Z}(\overline{\mathcal{A}})} \left( P_0(\mathcal{X}_0(T) = Z) \sum_{A \in Z} \frac{c(A)}{\beta(A)} \right).$$

PROOF. For any attribute $A \in \overline{\mathcal{A}}$, let $\lambda_A(T) = P_0(A \in \mathcal{X}_0(T))$ and let $q_A(T) = P_0(X_A$ performed by $T)$. Note that we have two useful expressions for $\lambda_A(T)$,

$$(8) \qquad \lambda_A(T) = \sum_{\substack{Z \in \mathcal{Z}(\overline{\mathcal{A}}) \\ Z \ni A}} P_0(\mathcal{X}_0(T) = Z)$$

and

$$(9) \qquad \begin{aligned} \lambda_A(T) &= P_0(A \in \mathcal{X}(T), X_A = 0) \\ &= P_0(A \in \mathcal{X}(T))P_0(X_A = 0) = q_A(T)\beta(A), \end{aligned}$$

where the second equality comes from the fact that the event that $A$ is performed by $T$ only depends on the values of tests for other attributes, and is thus independent of $X_A$ by the independence assumption.

Now recalling expression (3) we have

$$\begin{aligned} E_0[C(T)] &= \sum_{A \in \overline{\mathcal{A}}} c(A)q_A(T) \\ &= \sum_{A \in \overline{\mathcal{A}}} \frac{c(A)}{\beta(A)} q_A(T)\beta(A) \\ &= \sum_{A \in \overline{\mathcal{A}}} \frac{c(A)}{\beta(A)} \lambda_A(T) \\ &= \sum_{A \in \overline{\mathcal{A}}} \frac{c(A)}{\beta(A)} \sum_{\substack{Z \in \mathcal{Z}(\overline{\mathcal{A}}) \\ Z \ni A}} P_0(\mathcal{X}_0(T) = Z) \\ &= \sum_{Z \in \mathcal{Z}(\overline{\mathcal{A}})} \left( P_0(\mathcal{X}_0(T) = Z) \sum_{A \in Z} \frac{c(A)}{\beta(A)} \right). \qquad \square \end{aligned}$$

This lemma combines two straightforward observations. First, the cost "generated" by a specific attribute $A$ using strategy $T$ can be written as

$$(10) \qquad \begin{aligned} c(A)P_0(A \in \mathcal{X}(T)) &= \frac{c(A)}{\beta(A)} P_0(A \in \mathcal{X}(T))P_0(X_A = 0) \\ &= \frac{c(A)}{\beta(A)} P_0(A \in \mathcal{X}_0(T)). \end{aligned}$$

Second, the sum over attributes of the last expression can be reformulated as a sum over coverings (using the no-error property). Note in particular that (10) has the following interpretation: As far as average cost is concerned, it is equivalent to (i) pay the cost $c(A)$ every time test $X_A$ is performed, or (ii) pay the cost $c(A)/\beta(A)$ when $X_A$ is performed *and returns the answer* 0 but pay nothing otherwise.

Note also that the lemma implies that the average cost $E_0[C(T)]$ is therefore a *convex combination* of the quantities $\sum_{A \in Z} \frac{c(A)}{\beta(A)}$ for $Z \in \mathcal{Z}(\overline{\mathcal{A}})$.

5.2. *Filtering one special pattern.* Recall this corresponds to the case where the set of attributes has the structure of a vine (see Figure 2). We can imagine two broad scenarios: In one case, there is really only one pattern of interest, and hence no issue of invariance other than guaranteeing that every test is positive whenever $Y = y^*$. Imagine, for example, constructing a sequence of increasingly precise "templates" for a given shape, in which case *both* power *and* cost would typically increase with precision. In another scenario, one could imagine utilizing a hierarchy of tests originally constructed for multiple patterns in order to check for the presence of a single pattern $y^*$. Clearly, only one particular branch of the hierarchy is then relevant, namely the branch along which all the attributes contain $y^*$. Obviously, such tests would typically be less dedicated to $y^*$ than in the first scenario, except at the final level. In either case the natural framework is a *sequence of tests*, say $X_\ell$ for attributes $A_\ell$, with costs $c_\ell$ and powers $\beta_\ell$ for $\ell = 1, \ldots, L$, and the natural background measure is conditional on $Y \neq y^*$. Also, it is simpler here to consider the augmented hierarchy setting, so that we assume that there is a test at level $L + 1$ with $\beta_{L+1} = 1$, $c_{L+1} = c^*$.

The important quantity is the cost normalized by the power, $\{\frac{c_\ell}{\beta_\ell}\}$. Let $n(\ell)$, $\ell = 1, \ldots, L + 1$, denote the ordering of these ratios,

$$\tag{11} \frac{c_{n(1)}}{\beta_{n(1)}} \leq \frac{c_{n(2)}}{\beta_{n(2)}} \leq \cdots \leq \frac{c_{n(L+1)}}{\beta_{n(L+1)}}.$$

Since we are in the setting of the augmented hierarchy, there exists a distinguished index $\ell^*$ corresponding to the perfect test for which $c_{n(\ell^*)} = c^*$, $\beta_{n(\ell^*)} = 1$.

THEOREM 1. *The optimal strategy for detecting a single target pattern is to order the tests in accordance with* $(n(1), n(2), \ldots, n(\ell^*))$, *that is, perform* $X_{n(1)}$ *first, then* $X_{n(2)}$ *whenever* $X_{n(1)} = 1$, *and so on, and stop with* $X_{n(\ell^*)}$. *The tests* $X_{n(k)}$ *for* $k > \ell^*$ *are never performed.*

Note that the last test, $X_{n(\ell^*)}$, is the perfect one, and always returns the answer 0 under $P_0$. Reinterpreted in the original model, this would mean that if $X_{n(\ell^*-1)}$ is reached in the strategy and returns answer 1, then the testing procedure ends and the postprocessing stage is performed.

This theorem is a consequence of a straightforward recursion (proof omitted) applied to the following lemma.

LEMMA 2. *There exists an optimal strategy for which the first test performed is $X_{n(1)}$.*

PROOF. Let $T$ be some strategy performing the tests in the order $n'(1), n'(2), \ldots, n'(k^*)$ [for some $k^* \leq L+1$, with $n'(k^*) = n(\ell^*) = L+1$]. Assume $n'(1) \neq n(1)$ and consider strategy $T_0$ obtained by "switching" $X_{n(1)}$ to the first position, that is, performing $X_{n(1)}$ first, and then whenever $X_{n(1)} = 1$ continuing through strategy $T$ normally except if an index $i$ is encountered for which $n'(i) = n(1)$, in which case $X_{n(1)}$ is not performed again, but just skipped.

Compare the costs of $T$ and $T_0$ using (10) summed over attributes: clearly the mean cost of these strategies is a convex combination of the $(c_\ell/\beta_\ell)$, $\ell = 1, \ldots, L+1$, since $\sum_{\ell=1}^{L+1} P(A_\ell \in \mathcal{X}_0(T)) = 1$ in the single-pattern case. More explicitly,

$$P(A_k \in \mathcal{X}_0(T)) = \beta_k \prod_{\ell : n'(\ell) < n'(k)} (1 - \beta_\ell)$$

with the corresponding formula for $T_0$. From this formula it is clear that the weight for the ratio $c_{n(1)}/\beta_{n(1)}$ is higher in $T_0$ than in $T$, while all the other weights either are smaller or stay the same (depending whether the corresponding tests were placed before or after $X_{n(1)}$ in $T$). Since $c_{n(1)}/\beta_{n(1)}$ is the smallest of the ratios, the average cost of $T_0$ is lower than the cost of $T$.  □

5.3. *Filtering all patterns.* Our goal is to determine conditions under which (4) is minimized by the CTF strategy. First, we consider a simple sufficient condition which guarantees that the optimal strategy is complete, meaning $T \in \overline{\mathcal{T}}$. [Recall that $T \in \overline{\mathcal{T}}$ if $\widehat{Y}(T) = \widehat{Y}(\mathcal{X})$; in other words, testing is halted if and only if all "1-chains" in $\mathcal{X}$ are determined.] This condition is by no means necessary since we will prove the optimality of the CTF strategy (which belongs to $\overline{\mathcal{T}}$) under a much weaker condition, but is, however, informative.

PROPOSITION 1. *If for any attribute $A \in \mathcal{A}$, $\frac{c(A)}{\beta(A)} \leq c^*$, then the optimal strategy must belong to $\overline{\mathcal{T}}$.*

PROOF. Let $T$ be an optimal strategy and let $s$ denote a leaf of $T$. Recall that $\mathcal{X}(s)$ is the set of tests along the branch terminating in $s$ and $\widehat{Y}(s) = \mathcal{Y} \setminus \bigcup \{X_A \in \mathcal{X}(s) | X_A = 0\}$. The expected cost of $T$ is then of the form

$$(12) \qquad E_0[C(T)] = C + p_s c^* |\widehat{Y}(s)|,$$

where $p_s$ is the probability of reaching $s$, the second term is the contribution to the mean postprocessing cost at leaf $s$ and $C$ denotes the contributions of other nodes to the average cost. In general $\widehat{Y}(\mathcal{X}) \subset \widehat{Y}(s)$, and if these sets do not coincide, then by definition there must be a test $X_A \notin \mathcal{X}(s)$ for which $A \cap \widehat{Y}(s) \neq \varnothing$. Consider the strategy $T'$ obtained by adding this test to $T$ at node $s$. Then

$$(13) \qquad E_0[C(T')] = C + p_s[c(A) + \beta(A)c^*|\widehat{Y}(s) \setminus A| + (1 - \beta(A))c^*|\widehat{Y}(s)|].$$

Since $|\widehat{Y}(s)| - |\widehat{Y}(s) \setminus A| \geq 1$ it follows easily from the hypothesis, (12) and (13) that $E_0[C(T)] - E_0[C(T')] > 0$, which contradicts the optimality of $T$. □

We now turn to the problem of optimality of CTF strategies. The method of proof used in Section 5.1, although very simple in that case, will still serve as a template for most of the results to come. More precisely, under different assumptions about the models, we will always try to first establish the following property, denoted (CF) for "coarsest first":

DEFINITION 6 [(CF) *property*]. Test hierarchy $\mathcal{X}$ satisfies the (CF) property if there exists an optimal strategy for which the first test performed is the coarsest one.

In most cases, we will establish the optimality of $T_{\text{ctf}}$ as a consequence of (CF) for the various models considered. The current model—fixed, power-based cost—is the simplest and allows us to present the main ideas behind the arguments based on the (CF) property—a recursion based on "subhierarchies" and the concept of a "conditional strategy." As always, $\mathcal{A}$ is a nested hierarchy of attributes.

DEFINITION 7 (*Subhierarchy*). We call $\mathcal{B} \subset \mathcal{A}$ a subhierarchy of $\mathcal{A}$ if there exists an attribute $B_0 \in \mathcal{A}$ such that

$$\mathcal{B} = \{A \in \mathcal{A} | A \subseteq B_0\}.$$

More specifically, we call $\mathcal{B}$ the *subhierarchy rooted in* $B_0$ and we refer to $B_0$ as the set of patterns spanned by $\mathcal{B}$, also denoted $\mathcal{Y}_{\mathcal{B}}$.

DEFINITION 8 (*Conditional strategy*). Let $A_1$ be the root of $\mathcal{A}$ and let $\mathcal{B}$ be a subhierarchy of $\mathcal{A}$ rooted in one of the children of $A_1$. Then $\mathcal{A}$ can be written as a disjoint union $\mathcal{A} = \{A_1\} \dot{\cup} \mathcal{B} \dot{\cup} \overline{\mathcal{B}}$. Let $x_{\overline{\mathcal{B}}}$ be a set of numbers in $\{0, 1\}$ indexed by $\overline{\mathcal{B}}$. Consider a testing strategy $T$ for $\mathcal{A}$. The conditional strategy $T_{\mathcal{B}}(x_{\overline{\mathcal{B}}})$ on subhierarchy $\mathcal{B}$ is defined as follows: For every internal node $t$ of $T$:

(i) If $X(t)$ is a test for an attribute $B \in \mathcal{B}$, leave it unchanged.

(ii) If $X(t) = X_{A_1}$, cut the strategy subtree rooted at $t$ and replace it by the right subtree of $t$.

(iii) If $X(t)$ is a test for an attribute $A \in \overline{\mathcal{B}}$, cut the strategy subtree rooted at $t$ and replace it by the right subtree of $t$ if $x_A = 1$, and by the left subtree of $t$ if $x_A = 0$.

Finally, relabel every remaining leaf $s$ by $\widehat{Y}(s) \cap \mathcal{Y}_{\mathcal{B}}$.

This rather involved definition simply says that $T_{\mathcal{B}}(x_{\overline{\mathcal{B}}})$ is the testing strategy on subhierarchy $\mathcal{B}$ obtained from $T$ when $X_{A_1} = 1$ and the answers to $X_{\overline{\mathcal{B}}} = \{X_B, B \in \overline{\mathcal{B}}\}$ are fixed to be $x_{\overline{\mathcal{B}}}$, and $T$ is pruned accordingly. An obvious but nevertheless crucial observation is that $T_{\mathcal{B}}(x_{\overline{\mathcal{B}}})$ is indeed a valid testing strategy for the subset of attributes $\mathcal{B}$ and the corresponding subset of patterns $\mathcal{Y}_{\mathcal{B}}$.

THEOREM 2. *If property* (CF) *holds for any subhierarchy* $\mathcal{B}$ *of* $\mathcal{A}$ *(including* $\mathcal{A}$ *itself), then the CTF strategy is optimal.*

PROOF. The proof is based on a simple recursion. Let $L$ be the depth of $\mathcal{A}$. The case $L = 1$ is obvious from the (CF) property. Suppose the theorem is valid for any $L < L_0$ with $L_0 \geq 2$. Now consider the case $L = L_0$.

Let $T$ be an optimal testing strategy. From the (CF) property, we can assume that the test at the root of $T$ is $X_{A_1}$, the attribute at the root of $\mathcal{A}$. Denote by $\mathcal{B}_1, \ldots, \mathcal{B}_k$ the subhierarchies rooted at the children of $A_1$, which are of depth at most $L_0 - 1$. Since $\mathcal{A} = \{A_1\} \,\dot\cup\, \mathcal{B}_1 \,\dot\cup\, \cdots \,\dot\cup\, \mathcal{B}_k$ (a disjoint union), and $\mathcal{Y} = \mathcal{Y}_{\mathcal{B}_1} \,\dot\cup\, \cdots \,\dot\cup\, \mathcal{Y}_{\mathcal{B}_k}$, we can partition the cost of $T$ as follows:

$$
\begin{aligned}
E_0[C(T)] &= \sum_{A \in \mathcal{A}} q_A(T) c(A) + E_0[c^* | \widehat{Y}(T)|] \\
(14) \qquad &= q_{A_1}(T) c(A_1) + \sum_{A \in \mathcal{B}_1} q_A(T) c(A) + E_0[c^* | \widehat{Y}(T) \cap \mathcal{Y}_{\mathcal{B}_1}|] \\
&\quad + \cdots + \sum_{A \in \mathcal{B}_k} q_A(T) c(A) + E_0[c^* | \widehat{Y}(T) \cap \mathcal{Y}_{\mathcal{B}_k}|].
\end{aligned}
$$

Let us focus on the first sum. Let $\Omega(x_{\overline{\mathcal{B}_1}})$ be the event $\{X_{\overline{\mathcal{B}_1}} = x_{\overline{\mathcal{B}_1}}, X_{A_1} = 1\}$. Consider the conditional strategy $T^{(1)} = T_{\mathcal{B}_1}(x_{\overline{\mathcal{B}_1}})$ and let $q_A(T^{(1)}; x_{\overline{\mathcal{B}_1}})$ be the probability under $P_0(\cdot | \Omega(x_{\overline{\mathcal{B}_1}}))$ of performing the test for $A \in \mathcal{B}_1$ using $T^{(1)}$. The tests $\{X_A, A \in \mathcal{B}_1\}$ are conditionally independent given $\Omega(x_{\overline{\mathcal{B}_1}})$, with powers $\{\beta_A, A \in \mathcal{B}_1\}$. By the recurrence hypothesis, we can apply the theorem to subhierarchy $\mathcal{B}_1$ for the above conditional probability and conclude that the cost of strategy $T^{(1)}$ satisfies [for any $(x_{\overline{\mathcal{B}_1}})$]

$$
\begin{aligned}
E_0[C(T^{(1)}) | \Omega(x_{\overline{\mathcal{B}_1}})] &= \sum_{A \in \mathcal{B}_1} c(A) q_A(T^{(1)}; x_{\overline{\mathcal{B}_1}}) + E_0[c^* | \widehat{Y}(T^{(1)})| | \Omega(x_{\overline{\mathcal{B}_1}})] \\
(15) \qquad &\geq E_0[C(T_{\mathrm{ctf}}^{(1)}) | \Omega(x_{\overline{\mathcal{B}_1}})] = E_0[C(T_{\mathrm{ctf}}^{(1)})],
\end{aligned}
$$

where $T_{\mathrm{ctf}}^{(1)}$ is the CTF strategy for hierarchy $\mathcal{B}_1$ [whose cost is independent of $\Omega(x_{\overline{\mathcal{B}_1}})$]. We now want to take the expectation of (15) conditional on $\{X_{A_1} = 1\}$ only; by independence of the tests this amounts to taking the expectation of (15) with respect to $(X_{\overline{\mathcal{B}_1}})$. Now, by construction of the conditional strategy, denoting by $\beta_1$ the power of test $X_{A_1}$, for all $A \in \mathcal{B}_1$ we have

$$E_0[q_A(T^{(1)}; X_{\overline{\mathcal{B}_1}}) | X_{A_1} = 1] = P_0[X_A \text{ performed by } T | X_{A_1} = 1]$$
$$= q_A(T)(1 - \beta_1)^{-1},$$

where the last equality holds because $X_{A_1}$ is the first test to be performed in $T$. Similarly, on the event $\{X_{A_1} = 1\}$ we have $\widehat{Y}(T^{(1)}) = \widehat{Y}(T) \cap \mathcal{Y}_{\mathcal{B}_1}$, and therefore

$$E_0[c^* | \widehat{Y}(T^{(1)}) | | X_{A_1} = 1] = E_0[c^* | \widehat{Y}(T) \cap \mathcal{Y}_{B_1} | | X_{A_1} = 1]$$
$$= E_0[c^* | \widehat{Y}(T) \cap \mathcal{Y}_{B_1} |](1 - \beta_1)^{-1},$$

[since $X_{A_1} = 0 \Rightarrow \widehat{Y}(T) = \varnothing$]. Therefore, taking expectations w.r.t. $(X_{\overline{\mathcal{B}_1}})$ in (15) we obtain

$$E_0[C(T_{\mathcal{B}_1}(X_{\overline{\mathcal{B}_1}}))] = (1 - \beta_1)^{-1}\left( \sum_{A \in \mathcal{B}_1} c(A)q_A(T) + E_0[c^* | \widehat{Y}(T) \cap \mathcal{Y}_{B_1} |] \right)$$
$$\geq E_0[C(T_{\mathrm{ctf}}^{(1)})].$$

Applying the same reasoning to the other terms of (14), we now obtain

$$E_0[C(T)] \geq c(A_1)q_{A_1}(T) + (1 - \beta_1)E_0[C(T_{\mathrm{ctf}}^{(1)}) + \cdots + C(T_{\mathrm{ctf}}^{(k)})].$$

Finally, the right-hand side is precisely the total cost of the CTF strategy for $\mathcal{A}$. Therefore the CTF strategy is optimal. □

We now give a sufficient condition ensuring the (CF) property.

THEOREM 3. *Let $A_1$ be the coarsest test. Then the* (CF) *property holds under the condition*

$$\frac{c(A_1)}{\beta(A_1)} \leq \inf_{Z \in \mathcal{Z}(\overline{\mathcal{A}})} \sum_{A \in Z} \frac{c(A)}{\beta(A)}.$$

COROLLARY 1. *Consider the augmented hierarchy $\overline{\mathcal{A}}$ as a tree structure (the original hierarchy $\mathcal{A}$ can then be seen as the set of internal nodes of $\overline{\mathcal{A}}$). For any $A \in \mathcal{A}$, let $\mathcal{C}(A)$ be the set of direct children of $A$ in $\overline{\mathcal{A}}$. Then the CTF strategy is optimal if the following condition is satisfied:*

$$(16) \qquad \forall A \in \mathcal{A} \qquad \frac{c(A)}{\beta(A)} \leq \sum_{B \in \mathcal{C}(A)} \frac{c(B)}{\beta(B)}.$$

PROOF OF THEOREM 3. For this proof, it is easier to work with the "augmented" model put forward in Section 4.5.4. Let $T$ be a testing strategy for $\overline{\mathcal{A}}$ such that the first attribute to be tested is not the coarsest attribute $A_1$. From $T$ construct the strategy $T_0$ by "switching" test $X_{A_1}$ to the root, that is, perform $X_{A_1}$ first, and when the result is 1, proceed normally through strategy $T$, except when test $X_{A_1}$ is encountered in $T$, in which case it is not performed again and one jumps directly to its right child (corresponding to $X_{A_1} = 1$ in the original $T$).

Now compare the means cost of $T$ and $T_0$ using (7). Similarly to the proof of Lemma 2, we will prove that in the convex combination defining the cost in (7), the weight of the term $c(A_1)/\beta(A_1)$ is higher in $T_0$ than in $T$, while the weights of all the other terms of the form $(\sum_{A \in Z} c(A)/\beta(A))$ are smaller or stay unchanged for all other coverings $Z \in \mathcal{Z}(\overline{\mathcal{A}})$. This together with the hypothesis of the theorem establishes property (CF).

To verify the above statements about the weights of the different coverings, first call the "covering support" $\mathrm{CS}(T)$ of a strategy $T$ the set of coverings $Z \in \mathcal{Z}(\overline{\mathcal{A}})$ such that $P_0(\mathcal{X}_0(T) = Z) \neq 0$. It is clear from the construction of $T_0$ that $\mathrm{CS}(T_0) \subset \mathrm{CS}(T) \cup \{\{A_1\}\}$. Therefore we can restrict the analysis to the coverings in $\mathcal{Z}_0 = \mathrm{CS}(T) \cup \{\{A_1\}\}$.

Note that $\mathrm{CS}(T)$ is in one-to-one correspondence with the set of leaves of $T$ having nonzero probability to be reached; for any $Z \in \mathrm{CS}(T)$, $P(\mathcal{X}_0(T) = Z)$ is precisely the probability to reach the leaf $s_T(Z)$ of $T$ associated with the covering $Z$. Along the branch leading to this leaf one finds all the events $\{X_A = 0\}$ for $A \in Z$, along with a number of other events $\{X_A = 1\}$ for $A$ in a certain set $\mathcal{X}_1(s_T(Z))$. Therefore this probability is of the form

$$P_0(\mathcal{X}_0(T) = Z) = P_0(s_T(Z) \text{ is reached}) = \prod_{A \in Z} \beta(A) \prod_{A' \in \mathcal{X}_1(s_T(Z))} (1 - \beta(A')).$$

Now with this formula in mind, any $Z \in \mathcal{Z}_0$ falls into one of the following cases:

1. $Z = \{A_1\}$, in which case obviously $P_0(\mathcal{X}_0(T_0) = Z) \geq P_0(\mathcal{X}_0(T) = Z)$;
2. $A_1 \in Z$ but $Z \neq \{A_1\}$, in which case $P_0(\mathcal{X}_0(T_0) = Z) = 0$;
3. $A_1 \notin Z$ and $A_1 \notin \mathcal{X}_1(s_T(Z))$, in which case $P_0(\mathcal{X}_0(T_0) = Z) = (1 - \beta_1) \times P_0(\mathcal{X}_0(T) = Z)$;
4. $A_1 \notin Z$ and $A_1 \in \mathcal{X}_1(s_T(Z))$, in which case $P_0(\mathcal{X}_0(T_0) = Z) = P_0(\mathcal{X}_0(T) = Z)$.

Together, these different cases prove the desired property: $\{A_1\}$ is the only covering having higher weight in the cost of $T_0$ than in the cost of $T$. $\square$

Corollary 1 follows immediately: Its hypothesis clearly implies that the hypothesis of Theorem 3 is satisfied for any subhierarchy of $\mathcal{A}$ and the conclusion then follows from Theorem 2.

Note that, in contrast to what happened in the case of single-target detection, condition (16) falls short of being a necessary condition for ensuring the optimality of CTF strategies. To obtain a counterexample, consider the case of a depth-2 hierarchy with a coarsest attribute $A_1$ and two children $B_1$, $B_2$, and suppose that $c^*$ is large enough so that the condition of Proposition 1 is satisfied, so that we may restrict our attention to complete strategies. Then one can show (by explicitly listing all possible strategies) that the CTF strategy is optimal iff

$$\frac{c(A_1)}{\beta(A_1)} \le \inf\left(\frac{c(B_1)}{\beta(B_1)\beta(B_2)} + \frac{c(B_2)}{\beta(B_2)}, \frac{c(B_1)}{\beta(B_1)} + \frac{c(B_2)}{\beta(B_1)\beta(B_2)}\right).$$

Clearly this condition is weaker than (16).

*Application to the power-based cost model.* We can now look at the consequences of these results if we assume the cost model given by (5), in which case the following corollary is straightforward.

COROLLARY 2. *Assume the cost of the attribute tests obeys the model given by* (5), *with* $\Gamma$ *subadditive and* $\Psi(x)/x$ *increasing. Then the CTF strategy is optimal for any hierarchy* $\overline{\mathcal{A}}$ *for which* $\beta(A) \le \beta(B)$ *whenever* $B \subset A$. *In that case, the optimal strategy is CTF in both resolution and power.*

Similarly, in the case of detecting a single pattern of interest, if we assume $\Gamma \equiv 1$, the CTF strategy is optimal when $\Psi(x)/x$ is increasing, a result that was already proved in [13].

5.4. *Simulations with an elementary dependency model.* We also performed limited simulations in the case where the tests are not independent under $P_0$ but obey a very simple Markov dependency structure. Suppose the power of the coarsest test is $\beta_1$; the powers of subsequent tests follow a first-order Markov model depending on their direct ancestor. More precisely, the probability that a test returns 0 is $\gamma$ (resp. $\lambda$) given that its father returned 0 (resp. 1) with $\gamma \ge \lambda$. The cost model used is the multiplicative cost model given $c(X_{A,\beta})$ in (5), with $\beta$ the average power of the test.

We performed experiments for a set of four patterns and a corresponding depth-3 dyadic hierarchy, comparing the cost of the CTF strategy to the best cost among a set of 5000 randomly sampled strategies. In our experience, due to the restrained size of the problem, when there are in fact strategies better than the CTF one, then this is usually detected in the simulation.

What we found was that, for a given value of $\gamma$ and $\lambda$, the CTF strategy is generally optimal when $\beta_1 \le \lambda$ (for various choices of the power function $\Psi$). However, when $\beta_1$ becomes too large, then the CTF strategy is no longer optimal. Heuristically, this is because the coarse questions are then more powerful but also much too costly. The limiting value of $\beta_1$ for which CTF is optimal does not appear to be equal to the value $\beta^* = \lambda/(1 + \lambda - \gamma)$, the invariant probability for the Markov model. In particular, there are cases where $\lambda < \beta_1 \leq \beta^*$ (meaning that the average powers are increasing with depth) and yet CTF is not optimal.

To conclude, these very limited simulations seem to suggest that, even though the optimization problem is already somewhat complex even with a simple dependency structure and leads to challenging questions, still the optimality of CTF strategies can be expected to persist over a fairly wide range of models.

## 6. Optimal strategies for power-based cost and variable powers.

6.1. *Model and motivations.* In this section we only consider searching for all possible patterns. The previous section dealt with a fixed hierarchy—a single test $X_A$ at a given power $\beta(A)$ for each $A \in \mathcal{A}$. Now suppose we can have, for each $A \in \mathcal{A}$, tests of varying power; of course, a more powerful test at the same level of invariance will be more expensive. (In Section 8 we illustrate this trade-off for a particular data-driven construction.) In fact, for each attribute $A \in \mathcal{A}$, we suppose there is a test for every possible power, whose cost is determined as follows:

*Cost model.* Let $\Psi : [0,1] \to [0,1]$ be convex and strictly increasing with $\Psi(0) = 0$ and $\Psi(1) = 1$ and let $\Gamma : \mathbb{N}^* \to \mathbb{R}_+$ be subadditive with $\Gamma(1) = 1$. We suppose

$$(17) \qquad c(X_{A,\beta}) = c(A, \beta) = c \times \Gamma(|A|) \times \Psi(\beta).$$

Recall that the total cost of a strategy $T$ is given by

$$C_{\text{test}}(T) + c^* |\widehat{Y}(T)|.$$

The constant $c$ in (17) represents the cost of a $P_0$-perfect test for a single pattern and the constant $c^*$ represents the cost per pattern of disambiguating among the patterns remaining after detection. Evidently, only the ratio $c/c^*$ matters. We are going to assume that $c^* = c = \Psi(1) = 1$; note that this choice coheres with the formal interpretation of postprocessing cost as the cost of "errorless testing" put forward in Section 4.5.4.

For the rest of this section we will implicitly adopt this point of view, that is, replacing effective postprocessing cost by formal perfect tests corresponding to an additional layer of formal attributes copying the original leaves (this formal doubling of the leaf attributes allows us to keep untouched the rule that no attribute can be tested twice). For these special tests only, the power cannot be chosen arbitrarily and is fixed to 1; and the strategies considered must make no errors, enforced by performing at the end of the search some of these perfect tests if needed.

We are going to focus primarily on the case $\Gamma(k) = k$. Consequently,

$$(18) \qquad c(X_{A \cup B, \beta}) = c(X_{A,\beta}) + c(X_{B,\beta}) \qquad \text{when } A \cap B = \varnothing.$$

This is, in effect, the choice of $\Gamma$ least favorable to CTF strategies since there is no savings in cost due to shared properties among disjoint attributes. For instance, in practice it should not be twice as costly to build a test at power $\beta$ for the explanation $\{E, F\}$ as for $\{E\}$ or $\{F\}$ separately at power $\beta$, since (upon registration) these shapes share many "features" (e.g., edges; see Section 8). Nonetheless, with this choice of $\Gamma$ the convexity assumption for $\Psi$ can now be justified as follows:

*Motivation for convexity.* As usual, two tests for disjoint attributes are independent under $P_0$. Consider the following situation: For $A$ and $B$ disjoint, first test $A$ with power $\beta_1$ and stop if the answer is positive ($X_A = 1$); otherwise, test $B$ with power $\beta_2$ and stop. This produces a randomized, composite test for $A \cup B$ with power $\beta_1 \beta_2$ and (mean) cost

$$|A|\Psi(\beta_1) + \beta_1 |B| \Psi(\beta_2).$$

Contrast this with directly testing $A \cup B$ with power $\beta_1\beta_2$, which should not have greater cost than the composite test since, presumably, we have already selected the "best" tests at any given power and invariance; see Section 8 for an illustration. Under our cost model, this implies

$$(19) \qquad (|A| + |B|)\Psi(\beta_1\beta_2) \le |A|\Psi(\beta_1) + \beta_1 |B| \Psi(\beta_2).$$

Demanding (19) for any two attributes implies (by letting $|A|/|B| \to 0$) that we should have

$$(20) \qquad \Psi(\beta_1\beta_2) \le \beta_1 \Psi(\beta_2).$$

[Conversely, it is easy to see that if (20) is satisfied, then (19) holds for any $|A|, |B|$.] Since we want (20) to hold for any $\beta_1, \beta_2 \in [0, 1]$ we see (after dividing by $\beta_1\beta_2$) that (20) implies that $\Psi(x)/x$ is an increasing function. In our model we make the stronger hypothesis that $\Psi$ is convex in order to simplify the analysis.

*Remark on independence.* It would be unrealistic to assume the independence of *all* the tests in the variable-power hierarchy $\widetilde{\mathcal{X}}$, rather than for families corresponding to different attributes. In fact, in practice there is a limit to the number of independent (or even weakly dependent) tests that can be made for a fixed attribute. Were there not, then near-perfect detection would be possible in the sense of obtaining arbitrarily low cost and error by performing enough cheap tests of high invariance, at least in the case in which $\Psi'(0) = 0$.

EXAMPLE. Let $A = \mathcal{Y}$, the coarsest attribute, and suppose $\{X_{A,\beta_j}, j = 1, 2, \ldots\}$ are independent with $\beta_j \searrow \delta$. Consider the vine-structured testing strategy $T_n$ which successively executes $X_{A,\beta_j}, j = 1, 2, \ldots, n$, stopping [with label $\widehat{Y}(T_n) = \varnothing$] as soon as a null response is found and otherwise yielding $\widehat{Y}(T) = \mathcal{Y}$. Then it is easy to show (see [6]) that $P_0(\widehat{Y}(T_n) = \mathcal{Y}) \leq (1 - \delta)^n$ and that $E_0[C_{\text{test}}(T_n)] \leq |\mathcal{Y}| \frac{\Psi(\beta_1)}{\delta}$. Since $\Psi(\delta)/\delta \to 0$, given $\varepsilon > 0$, we can choose $n$, $\delta$ and $\beta_1$ close enough to $\delta$ such that $P_0(\widehat{Y}(T_n) \neq \varnothing) < \varepsilon$ and $E_0[C(T_n)] = E_0[C_{\text{test}}(T_n)] + E_0[|\widehat{Y}(T_n)|] < \varepsilon$.

6.2. *Basic results.* In the sequel, $\Psi^*$ will denote the Legendre transform of $\Psi$,

$$\Psi^*(x) = \sup_{\beta \in [0,1]} (x\beta - \Psi(\beta)).$$

In addition, for any $a > 0$, define

$$\Psi_a^*(x) = a\Psi^*\left(\frac{x}{a}\right),$$

$$\Phi_a(x) = x - \Psi_a^*(x).$$

6.2.1. *Optimal power selection.* Consider partially specifying a strategy $T$ by fixing the attribute $A$ to be tested at each (internal) node but not the power. What assignment of powers (to the nonperfect) tests minimizes the average cost of $T$? As with dynamic programming, it is easily seen that the answer is given as follows: Start by optimizing the powers of the last, nonperfect tests performed along each branch (since the left and right subtrees of such a node have fixed, known cost), and then climb recursively up each branch of the tree, optimizing the power of the parent at each step. The actual optimization at each step is a simple calculation, summed up by the following lemma:

LEMMA 3. *Consider a (sub)strategy $T$ consisting of a test $X_{A,\beta}$ at the root, a left subtree $T_L$ of average cost $x$ and a right subtree $T_R$ of average cost $y$. Let $\Gamma(|A|) = a$. Then under the cost model* (17) *the average cost of $T$ using the optimal choice of $\beta$ is given by*

$$(21) \qquad E_0[C(T)] = y - \Psi_a^*(y - x) = x + \Phi_a(y - x).$$

*In particular, if $T_L$ is empty, then $x = 0$ and $E_0[C(T)] = \Phi_a(y)$. If $\Psi$ is differentiable, the optimal choice of $\beta$ is*

$$\beta^* = \begin{cases} (\Psi')^{-1}((y-x)/a), & \text{if } (y-x)/a \in \Psi'([0,1]), \\ 0, & \text{if } (y-x)/a < \Psi'(0), \\ 1, & \text{if } (y-x)/a > \Psi'(1). \end{cases}$$

*More generally, $\Psi$ admits $(y-x)/a$ as a subgradient at point $\beta^*$.*

PROOF. Let $T(\beta)$ denote the strategy using power $\beta$, and calculate the average cost of $T(\beta)$ as a function of $\beta, x, y, a$,

$$E_0[C(T(\beta))] = c(X_{A,\beta}) + \beta x + (1 - \beta)y$$
$$= a\Psi(\beta) + \beta(x - y) + y$$
$$= y - a\left(\left(\frac{y - x}{a}\right)\beta - \Psi(\beta)\right).$$

Now minimizing over $\beta$ leads directly to (21) and the formulae for $\beta^*$, using the definitions of $\Psi^*$ and $\Psi_a^*$. $\quad\square$

6.2.2. *Properties of the CTF strategy.* In previous sections, with fixed powers, all variations on CTF exploration (e.g., depth-first and breadth-first) had the same average cost, and hence we spoke of "the" CTF strategy. With variable powers the situation might appear different: The bottom-up optimization process in Section 6.2.1 for assigning the powers may lead to different mean costs for different CTF strategies. More specifically, recall that $A(s), \beta(s)$ denote the attribute and power assigned to an internal node $s$ in a tree $T$. For CTF trees, it may be that $\beta^*(s)$, the optimal power at $s$, depends on the position of $s$ within $T$ as well as $A(s)$.

The following theorem states that, in fact, as in the fixed-powers case, among CTF strategies, the order of testing is irrelevant *when the powers are optimally chosen.* More precisely, the optimal power of a test depends only on the attribute being tested, specifically on the structure of the subhierarchy rooted at the attribute. Consequently, in CTF strategies a given attribute will always be tested at the same power, which means that CTF designs can be implemented by constructing only one test per attribute—a considerable practical advantage.

THEOREM 4. *For any CTF strategy $T$, and for any two nodes $s, t$ in $T$ with $A(s) = A(t)$, the optimal choices of powers are identical: $\beta^*(s) = \beta^*(t)$. In fact, the unique power assigned to an attribute $A \in \mathcal{A}$ depends only on the structure of subhierarchy $\mathcal{B}(A)$ rooted in $\mathcal{A}$. As a consequence, all CTF strategies have the same average cost.*

Whereas the principle of the proof is simple (a recursion on the size of $\mathcal{A}$), it does require some auxiliary notation, and hence we postpone it to the Appendix.

Turning to the cost of the CTF strategy, it can easily be computed recursively for regular attribute hierarchies and the simple complexity function $\Gamma(k) = k$. More precisely, we have the following theorem for dyadic hierarchies, in which $\beta_\ell^*(L)$ denotes the optimal power for the $2^{\ell-1}$ attributes at level $\ell = 1, \ldots, L$ for a hierarchy of total depth $L$.

THEOREM 5. *Let $C_L$ denote the average CTF cost of a regular, complete dyadic hierarchy of depth $L$. Then*

$$(22) \qquad C_{L+1} = \Phi_{2^L}\left(2C_L\right)$$

*with (formally) $C_0 = \Psi(1)/2$. Furthermore,*

$$C_L/2^{L-1} \searrow \Psi'(0), \qquad L \to \infty,$$

*and*

$$(23) \qquad \beta_1^*(L) \searrow 0, \qquad L \to \infty.$$

*Finally,*

$$(24) \qquad \beta_\ell^*(L) = \beta_1^*(L - \ell + 1), \qquad \ell = 1, \ldots, L,$$

*from which it follows that the CTF strategy is CTF in power, that is, power increases with depth.*

PROOF. Consider a (complete, dyadic) hierarchy of depth $L + 1$. The coarsest attribute has cardinality $|A_1| = 2^L$ and the (optimized, breadth-first) CTF strategy starts with the corresponding test. If $X_{A_1} = 0$, the search is over; if not, it is necessary to pay the mean cost for the two subhierarchies of depth $L$. We thus apply (21) with $x = 0$, $y = 2C_L$ to obtain (22). When $L = 1$ (one pattern), it is easy to check that we retrieve the right value of $C_1$ from (22) with $C_0 = \Psi(1)/2$ by noting that, in this case, $y = \Psi(1)$, which is the cost of a perfect test.

Let $U_L = C_L/2^{L-1}$. Then (22) can be rewritten as

$$U_{L+1} = \Phi_1(U_L),$$

which allows us to study the asymptotic behavior of $U_L$ when $L$ is large based on the function $\Phi_1(x) = x - \Psi^*(x)$. Since $\Psi$ is convex, it follows that $\frac{\Psi(\beta)}{\beta}$ is increasing, and hence $x\beta - \Psi(\beta) \le 0$ for all $0 \le \beta \le 1$ (with equality at $\beta = 0$) whenever $0 \le x \le \Psi'(0)$. Consequently, $\Psi^*(x) = 0$ and $\Phi_1(x) = x$ for $x \in [0, \Psi'(0)]$. Similarly, $\Phi_1(x) < x$ for $x > \Psi'(0)$. We have $U_0 = \Psi(1) \ge \Psi'(0)$ because $\Psi$ is convex, and hence since $\Phi_1$ is concave, $U_L \searrow \Psi'(0)$ as $L \to \infty$. Finally, from Lemma 3 we can also conclude that $\beta_1^*(L) = (\Psi')^{-1}(U_L \wedge \Psi'(1))$. The last assertion (24) of the theorem follows directly from Theorem 4. □

REMARK 1. We deduce from the above results that if $\Psi'(0) = \delta > 0$, we have $C_L \sim \delta 2^{L-1}$. If, on the other hand, $\Psi'(0) = 0$, then $C_L = o(2^{L-1})$. This should be compared to the strategy of performing only (all) the perfect tests, which costs $2^{L-1}$.

REMARK 2. Since the optimal powers are increasing with depth, if we now consider them as fixed we are in the framework of Corollary 2 ensuring that, for these choices of powers, the CTF strategy is indeed optimal.

REMARK 3. Note that the cost of individual tests (with optimal powers) may not vary monotonically with their depth; however, the *cumulated* cost of all tests at a given depth increases with depth.

6.3. *Is the CTF strategy optimal*? We have not been able to prove the optimality of the CTF strategy under general conditions on $\Psi$, but rather only for one specific example. This is disappointing because the simulations presented later in this section strongly indicate a more general phenomenon.

If we try to follow our usual method for proving optimality, it turns out that the most difficult step is actually to prove the (CF) property. Under the (CF) property, the optimality of CTF would readily follow—it suffices to follow the lines of the proof of Theorem 2 with minor adaptations, mainly replacing families $(X_A)_{A \in \mathcal{B}}$ by $(X_{A,\beta})_{A \in \mathcal{B}, \beta \in [0,1]}$.

One way to prove the (CF) property is to proceed iteratively, repeatedly applying the "switching property":

DEFINITION 9 (*Switching property*). A power function $\Psi$ has the switching property if any (sub)tree $T$ of the form shown on the left-hand side of Figure 4, with any powers, has a larger mean cost than the tree obtained by switching the two first tests of $T$ (shown on the right-hand side of Figure 4), with optimal powers. Using Lemma 3, this inequality may be expressed as follows:

$$\forall y \geq x \geq 0, \forall a \geq b \geq 0$$

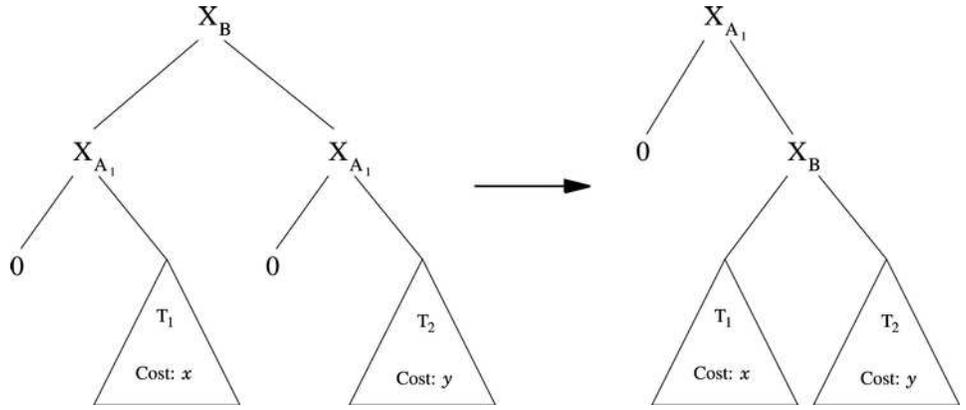

FIG. 4. *The context of the switching property. Attribute $A_1$ is the coarsest attribute in the hierarchy; hence $\Gamma(|B|) = b < \Gamma(|A_1|) = a$.*

(25)
$$\Phi_a(x + \Phi_b(y - x)) \le \Phi_a(x) + \Phi_b(\Phi_a(y) - \Phi_a(x)).$$

We then have the following lemma.

LEMMA 4.  *The* (CF) *property is implied by the switching property.*

PROOF.  Note first that we can assume that $X_{A_1}$, the coarsest test, is performed at some point (at some power) along every branch of any $T$. If this is not the case, it can simply be added, with zero power, at the end of any branch where it does not appear without changing the cost. Now let $T$ be a strategy such that $X_{A_1}$ is not performed *first*. Apply the switching lemma to any subtree of $T$ of the form shown on the left-hand side of Figure 4. In this way, $X_{A_1}$ is pushed up in the tree while reducing the cost. This can be done repeatedly until no such subtree exists, that is, the situation depicted in Figure 4 does not occur anywhere in $T$. But then the resulting tree must have $X_{A_1}$ at the root. Otherwise, let $k$ be the maximum depth in $T$ where $X_{A_1}$ appears, and let $s$ be the corresponding node. Let $s'$ be the direct sibling of $s$, which exists since $k > 1$. Consider a branch $b$ containing $s'$. Since $X_{A_1}$ is performed along any branch, it must be performed somewhere in $b$, say at node $t$. But $t$ cannot be an ancestor of $s'$, since otherwise $X_{A_1}$ would be performed twice along branch $b$, a contradiction. Nor can $t$ be a descendant of $s'$, since that would contradict the definition of $k$. Therefore $X_{A_1}$ is performed at $s'$, which contradicts the assumption that there is no subtree of the form shown on the left of Figure 4. This concludes the proof. □

From numerical experiments, we know, however, that the switching property is not satisfied for an arbitrary (convex) power function $\Psi$. Whereas we believe that it should be possible to prove the switching lemma under some additional conditions on $\Psi$, we have so far only been able to prove it for one case we refer to as the "harmonic" cost function,

(26)
$$\Psi(x) = 2 - 2\sqrt{1 - x} - x,$$

which we now investigate.

6.4. *CTF optimality for the harmonic cost function.*  Throughout this section $\Psi$ is given by (26). This function has the following properties:

1. $\Psi$ is convex and increasing;
2. $\Psi(0) = \Psi'(0) = 0$ and $\Psi(1) = 1$, $\Psi'(1) = \infty$;
3. $\Psi^*(x) = x - \frac{x}{x+1}$; $\Phi_a(x) = \frac{ax}{x+a} = (x^{-1} + a^{-1})^{-1}$.

Note that $x$ and $a$ have symmetric roles in $\Phi_a$, and that $\Phi_a(x)$ is the "harmonic sum" of $x$ and $a$.

We first study the switching lemma in the case of an empty left subtree $T_L$.

LEMMA 5. *Consider two tests $X_A$ and $X_B$ with $\Gamma(|A|) = a$ and $\Gamma(|B|) = b$. Let $T_{AB}$ be the tree shown on the right-hand side of Figure 4 with $T_1 = \varnothing$ and let $T_{BA}$ have the same structure with $X_A$ and $X_B$ reversed. Then, with the optimal assignment of powers to $X_A$ and $X_B$, both $T_{AB}$ and $T_{BA}$ have the same cost.*

PROOF. By applying Lemma 3 (with $x = 0$) twice, the cost of $T_{AB}$ is $\Phi_a \circ \Phi_b(y)$ and the cost of $T_{BA}$ is $\Phi_b \circ \Phi_a(y)$. It is then easy to check that

$$\Phi_a \circ \Phi_b(y) = \Phi_b \circ \Phi_a(y) = \frac{aby}{ay + by + ab} = (a^{-1} + b^{-1} + y^{-1})^{-1}. \qquad \square$$

NOTE. Clearly, $\Phi_a \circ \Phi_b(x)$ is the harmonic sum of $x, a$ and $b$. More generally, consider any "right vine" $T$ consisting of at most one test per level of resolution. Then, under $\Psi$ the average cost of $T$ (with optimal powers) is *independent* of the order in which the tests are performed; moreover, this average cost is simply the harmonic mean of the values $\Gamma(|A_i|)$ for the tests performed. In particular, this result is totally *independent* of the choice of the complexity function $\Gamma$.

We now return to the "full" switching lemma:

THEOREM 6. *The switching property—and hence the optimality of the CTF strategy—holds for the harmonic power function with any complexity function $\Gamma$.*

For the proof see the Appendix.

*Analogy with resistor networks.* We conclude this section with a curious connection: Consider a hierarchy of depth $L$ with coarsest attribute $A_1$ and $a_1 = \Gamma(|A_1|)$. Let $C_1$ be the average cost of the CTF strategy for the hierarchy with $A_1$ removed. From Lemma 3, with $x = 0$ and $y = C_1$,

$$E_0[C(T_{\text{ctf}})] = \Phi_{a_1}(C_1) = \frac{a_1 C_1}{a_1 + C_1} = \left(\frac{1}{a_1} + \frac{1}{C_1}\right)^{-1}.$$

This is exactly the conductance of an electrical circuit composed of two serial resistors of conductances $C_1$ and $a_1$. Continuing, $C_1$ is the sum of the CTF costs over the two subhierarchies of depth $L - 1$; if $C_1'$ denotes the cost of these hierarchies, the cost $C_1$ can be interpreted as the conductance of an electrical circuit formed from two parallel resistors, each of conductance $C_1'$. The global cost of the CTF strategy is therefore equal to the conductance of the tree-structured resistor network depicted in Figure 5 (wherein a row of resistors is added at the bottom of the tree in order to represent the cost of the postprocessing, or, equivalently, perfect testing). We observe that nothing would be changed in the case of a nonsymmetric, tree-structured hierarchy, even with attributes of varying complexities at the same level.

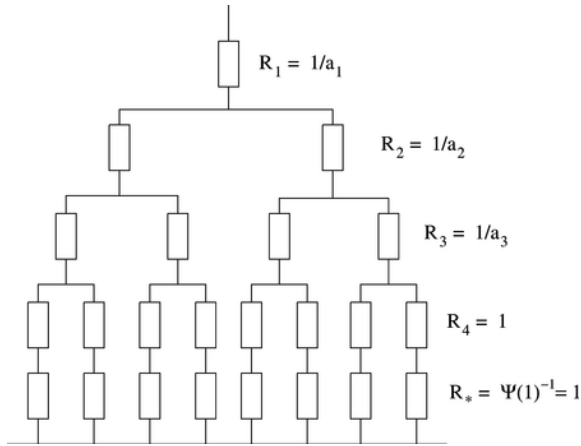

Fig. 5. *Tree-structured resistor network identified with the attribute hierarchy, where $a_l = \Gamma(|A_l|)$ is the complexity of attributes of level $l$ and $R_l = 1/a_l$ is the associated resistance; note that $a_4 = 1$ by convention for the bottom attributes. The last row of resistors represents the postprocessing stage. The conductance of this circuit is exactly the CTF testing cost of the attribute hierarchy when $\Psi$ is the harmonic power function.*

6.5. *Simulations.* In this section we investigate the optimality of CTF search by way of simulations involving several different power functions $\Psi$. In every case we take $\Gamma(k) = k$. The various choices of $\Psi$, and corresponding functions $\Phi_1(x) = x - \Psi^*(x)$, are presented in Table 1; obviously we have chosen functions with closed-form Legendre transforms. We took $\lambda = 1$ for $\Psi_4$ and $\mu = 8$ for $\Psi_7$.

First we investigated the switching property, which we know to be sufficient for the optimality of $T_{\text{ctf}}$. To this end, we computed and plotted the difference $\Delta(a, b, x, y)$ between the left-hand side and the right-hand side of the key inequality (25). Without loss of generality, we put $a = 1$. Plots of $\Delta(1, b, x, y)$ are given in [6] for the particular choice $b = 2$. The switching property is satisfied if the surface lies below the *xy-plane*. Some of these surfaces (corresponding to $\Psi_2, \Psi_4, \Psi_6$) clearly do not, whereas the others appear to satisfy this inequality (at least all sampled values are negative). In other experiments with other values of $b$ for $\Psi_1$ and $\Psi_3$ we always found

TABLE 1
*Convex power functions used in our simulations*

| Number | $\Psi$ | $\Phi_1$ |
|---|---|---|
| 1 | $x(1 - \sqrt{1-x})$ | $x - (1 - (1 - x + \frac{1}{9}\sqrt{(1-x)^2 + 3})^2)(\frac{2}{3}(x-1) + \frac{1}{3}\sqrt{(1-x)^2 + 3})$ |
| 2 | $x^2/2$ | $\begin{cases} x - x^2/2, & \text{if } x < 1, \\ \frac{1}{2}, & \text{otherwise} \end{cases}$ |
| 3 | $1 - \sqrt{1-x^2}$ | $1 + x - \sqrt{x^2 + 1}$ |
| 4 | $\exp(\lambda x) - 1$ | $\begin{cases} x, & \text{if } x < \lambda, \\ x - 1 - \frac{x}{\lambda}(\log(\frac{x}{\lambda}) - 1), & \text{if } \lambda \leq x \leq \lambda e^{\lambda}, \\ x - e^{\lambda} + 1, & \text{if } x > \lambda e^{\lambda} \end{cases}$ |
| 5 | $2 - x - 2\sqrt{1-x}$ | $x/(1+x)$ |
| 6 | $1 - \sqrt{1-x}$ | $\begin{cases} x, & \text{if } x < \frac{1}{2}, \\ 1 - \frac{1}{4x}, & \text{otherwise} \end{cases}$ |
| 7 | $\exp(\mu x) - 1 - \mu x$ | $\begin{cases} x(1 + \frac{1}{\mu}) - (1 + \frac{x}{\mu})\log(1 + \frac{x}{\mu}), & \text{if } x < \mu(e^{\mu} - 1), \\ e^{\mu} - 1 - \mu, & \text{otherwise} \end{cases}$ |

Note that $\Psi_5$ is the harmonic function.

$\Delta \leq 0$. However, we found regions with $\Delta > 0$ for $\Psi_7$ for higher values of $b$, and hence this cost function does not satisfy the switching property.

From these plots it is tempting to speculate that only power functions $\Psi$ such that $\Psi'(0) = 0$ and $\Psi'(1) = +\infty$ can satisfy the full switching property; however, these conditions are very likely not sufficient. Note that $\Psi'(0) = 0$ means that, at any given level of invariance, one can have an arbitrarily small cost-to-power ratio and $\Psi'(1) = +\infty$ means that very high powers are likely not worth the increased cost. Intuitively, both of these properties favor CTF strategies.

The second type of simulation was more direct. Strategies were sampled at random by the simplest method possible: we sampled purely attribute-based strategies $T$ by recursively visiting nodes and choosing an attribute $A \in \mathcal{A}$ at random subject to the two obvious constraints: (i) no attribute is repeated along the same branch, and (ii) no "useless" attribute is chosen, meaning that $A$ consists entirely of patterns already ruled out by the previous tests. Then, for each such $T$, powers were individually assigned to the tests at each node in order to minimize the cost, which was compared with that of the CTF strategy. This procedure was repeated for various choices of $\Psi$ [with $\Gamma(k) = k$] for regular, dyadic hierarchies for $|\mathcal{Y}| = 4$ patterns (i.e., $L = 3$) and for $|\mathcal{Y}| = 8$ patterns (i.e., $L = 4$). For each $\Psi$, we sampled several tens of thousands of trees $T$. [Of course the sheer number of possible strategies (modulo power assignments) in the case $L = 4$ is several orders of magnitude larger.] Summarizing our observations:

(a) *In all cases, the CTF strategy had lower cost than any other strategy sampled.*

(b) *Upon visual inspection, the best sampled strategies seemed close to the CTF strategy in the sense of only differing at relatively deep nodes.*

In conclusion, and bearing in mind the limited scope of both types of simulations, we believe the following conclusions are reasonable:

1. *The switching property is quite likely valid for cost models other than the harmonic function; however, it requires hypotheses in addition to convexity.*
2. *The optimality of the CTF strategy probably holds for a very wide range of cost models, including those which do not satisfy the switching property (for all values of $a, b, x, y$). As a result, requiring the switching property is likely too restrictive and, more generally, arguments based on the* (CF) *property may not be the most efficacious.*

**7. Remarks on a usage-based cost model.** In this section we summarize some results obtained in [6] for a somewhat different scenario. We consider only the case of a fixed-powers hierarchy. In this model, the cost of a test $c(X)$ may be chosen in accordance with the strategy employed; it depends on the "resource" $r(X)$ allocated to it [through a negative exponential function $r(X) = \exp(-c(X))$] and there is a global resource constraint, $\sum_{X \in \mathcal{X}} r(X) \le R \le 1$. This corresponds to the belief that in some circumstances it might not be efficient to fix the costs of the tests in advance, regardless of their inherent complexity. It may be more efficient to allow the utilization of computing resources to be partitioned in accordance with the frequency with which certain routines are performed; in this case the cost represents the computing time rather than the computing complexity. In this framework the optimal resource allocation gives rise to a *usage-based* cost; the cheapest tests are the ones used the most often in a given strategy. The testing cost of a strategy with optimal resource allocation is then (from standard arguments)

$$(27) \qquad E_0[C_{\text{test}}(T)] = -\sum_X q_X(T) \log(q_X(T)) + Q(T) \log(Q(T)/R),$$

where $Q(T) = \sum_X q_X(T)$. Furthermore, no postprocessing cost is taken into account, but we only allow complete strategies, so that the goal is to minimize (27) over complete strategies. In [6] we prove that for a hierarchy for which each attribute has at least two children and for which the powers are increasing with the resolution level, the CTF strategy is optimal if we assume that all tests have power greater than some constant $\beta_1 = 7/8$. While this (probably improvable) value is not entirely realistic as far as practical applications are concerned, we believe it is an important step in favor of CTF optimality for this cost model.

In addition, we argue that it makes sense in this latter framework to consider an extended scenario where repeated search tasks are undertaken for different sets of target patterns, whereas the resources are distributed in advance among all tests. While the set of targets changes from task to task, the individual attribute tests are reusable. The patterns are identified with conjunctions of abstract attributes at different resolution levels, taken from a possibly very large pool. Whereas the analysis in the fixed-cost model remains unchanged, there is a significant difference under usage-based cost since we must distribute the resources over a larger number of tests. In order to simplify the analysis, we suppose the set of target patterns $\mathcal{Y}$ is randomized for each new search task and again present some fairly mild sufficient conditions (about the dependence of power on resolution and the size of the attribute pools, e.g., exponential growth of pool size with resolution, and negative polynomial decrease of type II error) ensuring the optimality of CTF strategies.

**8. Applications to pattern recognition.**   In order to illustrate our framework for pattern recognition we present two types of results: First, we give a few examples of the scene interpretation problem and cite some previous work on a CTF strategy for object detection. Only pictures and references are provided. The purpose is merely to demonstrate the efficacy of the approach in a real computational vision problem. Second, in order to illustrate numerically the quantities appearing in our analysis, and to check whether the cost model is reasonable in at least one concrete setting, we outline a more or less exact implementation, due to Franck Jung, of the pattern filtering design for a synthetic example introduced in [20]—detecting rectangles amidst clutter. It was developed in order to automate cartography by detecting roofs of buildings in aerial photographs [21]. Only those aspects which shed light on the mathematical analysis are described; all the details may be found in [6].

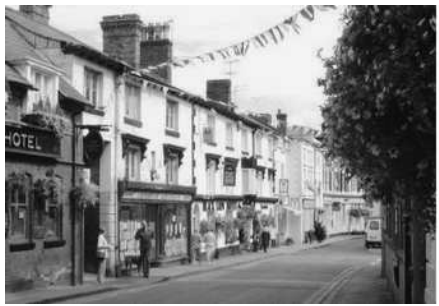 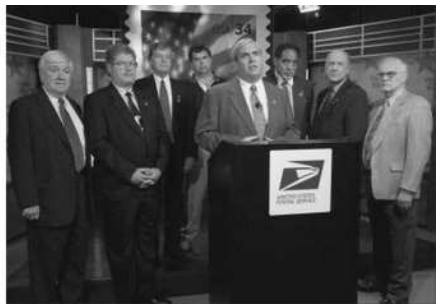

FIG. 6.   Left: *a "natural" image.* Right: *group photograph used in an experiment on face detection.*

8.1. *Scene interpretation.* Consider the scenes in Figure 6. The semantic interpretation of the left image (town, shops, pedestrians, etc.) is effortless for humans but far beyond what any artificial system can do. For the image on the right, the goal might be more modest—detect and localize the faces. Enriching the description with information about the precise pose (scale, orientation, etc.), identities or expressions would be more ambitious. Many methods have been proposed for face detection, including artificial neural networks [24], Gaussian models [26], support vector machines [22], Bayesian inference [9] and deformable templates [30].

To relate these tasks to the framework of this paper, imagine attempting to characterize a (randomly selected) subimage containing at most one object from a predetermined repertoire. (The whole scene can then be searched by a divide-and-conquer strategy; see Section 8.2 and [12].) The dominating explanation $Y = 0$ corresponds to "background" or "clutter" and each of the others, $Y \in \mathcal{Y}$, corresponds to the instantiation of an object wholly visible in the subimage. Even with only one (generic) object class, the number of possible instantiations is very large; that is, there is still considerable within-class variability. For instance, detecting a face at a fixed position, scale and orientation might not be terribly difficult, even given variations in lighting and nonlinear variations due to expressions; it can be accomplished with standard learning algorithms such as multilayer perceptrons, decision trees and support vector machines. However, the amount of computation required to do this separately for every possible pose is prohibitive. Instead, we propose to search simultaneously for many instantiations, say over a range of locations, scales and orientations. In our simplified mathematical analysis, that range of poses is $A = \mathcal{Y}$, which is the "scope" of our coarsest test $X_A$. (It may not be practical to envision a totally invariant test, in which case there are multiple hierarchies.)

This approach to scene interpretation has been shown to be highly effective in practice. A version involving successive partitions of object/pose pairings, rank-based tests for the corresponding (classes of) hypotheses and breadth-first CTF search appears in [19]. The detection results shown in Figure 7 were obtained by an algorithm [14] based on the strategy proposed here—traversing a multiresolution hierarchy of $\mathcal{X}$ binary hypothesis tests $\{X_A, A \in \mathcal{A}\}$, where each $A$ represents a family of shapes with some common properties and $X_A$ is an image functional designed to detect shapes in this family. In the face detection experiments, $A$ is a subset of affine poses and $X_A$ is based on checking for special local features (e.g., edges) which are likely to be present for faces with poses in $A$. In fact, $X_A$ can be interpreted as a likelihood ratio test [3]. Recently, researchers in the computer vision community have started using similar methods for similar problems; see, for example, [25] and [29]. Ideas related to CTF processing have also been proposed by [15] in a Bayesian classification framework where a hierarchy of estimators is built for the posterior of recursively clustered classes. In Figure 7, the efficiency of sequential testing is illustrated for the group photo by counting, for each pixel, the amount of computation performed in its vicinity; clearly the spatial "density of work" is highly skewed. The corresponding density would be flat for nearly all other methods, that is, those based on multilayer perceptrons or support vector machines.

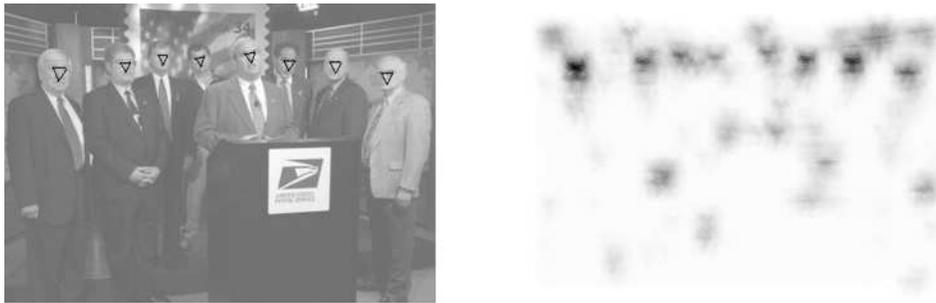

FIG. 7. *The detections (*left*) and "density of work" (*right*) for the group photo.*

8.2. *Rectangle detection.* The goal is to find and localize rectangles in a "scene" of the type shown in Figure 11. The generative model (which involves first inserting and degrading rectangles and then adding clutter) is described in [6].

There are many ways to find the rectangles. For instance, one could use any of the methods cited above for finding faces. For the artificial problem illustrated in Figure 11, with limited noise and clutter, it would not be surprising to obtain a decent solution with standard model-based or learning-based methods. Our intention is only to demonstrate how this might be done in an especially efficient manner with a sequential testing design.

8.2.1. *Problem formulation.* It is clearly impossible to find common but localized attributes of two rectangles with significantly different (geometric) poses, say far apart in the scene. Here, the "pose" of a rectangle has four parameters: orientation, center, height and length. Consequently, we divide the whole scene into nonoverlapping $5 \times 5$ regions and apply a simple, "divide-and-conquer" strategy based on location. Each $5 \times 5$ region $R$ is visited in order to determine if there is a rectangle in the scene whose distinguished point (say the center) lies in $R$; depending on its scale, the rectangle itself will enclose some portion of the scene surrounding $R$. We can assume that the scale of the rectangle is restricted to a given range whose lower end represents the smallest rectangles we attempt to find. Larger rectangles are found by repeatedly downsampling the image and parsing the scene in the same way; this is how the faces in Figure 6 were detected. (Similarly, the orientation of the rectangle is restricted to a given range of angles; other orientations could be found by repeating the process with suitably transformed detectors.)

Partitioning the scene into nonoverlapping regions and downsampling to handle scale can be thought of as the first two levels of a recursive partitioning of the full pose space. The loops over regions $R$ and scales are the "parallel component" of the algorithm and not of interest here. The serial component is a CTF search to determine if there is a rectangle within a range of scales whose center lies in a fixed region $R$. This is the heart of the algorithm and the real source of efficient computation. The hypothesis $Y = 0$ stands for "no rectangle with these parameters" and is evidently a complex mixture of configurations due to clutter, larger rectangles and nearby ones.

8.2.2. *Patterns, attributes and tests.* In order to define the set of explanations $\mathcal{Y}$, we partition the (reference) pose space $\Theta$ into small subsets. A "pattern" or "explanation" $y \in \mathcal{Y}$ is then a subset of poses at approximately the resolution of the pixel lattice. In fact, these subsets are, by definition, the cells at the finest layer of the attribute hierarchy—a recursive partitioning of $\Theta$ of the type used throughout the paper, yielding $\Theta = \{A_{l,k}\}$. In this case $Y$ represents the true pose at the pixel resolution.

There are $L = 6$ levels which correspond to five splits: two (binary) on orientation, one (quaternary) on position and two (binary) on scale (one on height and one on length). In particular there are $|\mathcal{Y}| = 64$ finest cells, each with resolution 1.25 pixels in location, two pixels in length and height, and $\pi/16$ radians in tilt. Let $\eta_l$ be the cardinality (scope) of attributes at resolution level $l$. The quaternary split happens to be the second one, and hence $(\eta_1, \ldots, \eta_6) = (64, 32, 8, 4, 2, 1)$.

As in the references cited above, the tests $X_A$ are extremely simple image functionals based on local features $\xi$ related to edges. Each test $X_A$ is based on a threshold $\tau = \tau(A)$ and a collection $\mathcal{S}(A)$ of these features (corresponding to varying positions, orientations and levels of resolution):

$$X_A = \begin{cases} 1, & \text{if } \sum_{\xi \in \mathcal{S}} \xi \geq \tau, \\ 0, & \text{otherwise.} \end{cases}$$

Thus, evaluating $X_A$ consists of checking for at least $\tau$ features among a special ensemble dedicated to $A$. Actually, we build many tests of varying powers for each $A \in \mathcal{A}$, each one corresponding to a different collection $\mathcal{S}$. Identifying $\mathcal{S}$ and $\tau$ is a problem in statistical learning. We use a fairly simple procedure which is described in [6].

The cost $c(X_A)$ is defined as the number of pixels involved in evaluating $X_A$, which is the number of pixels which participate in the definition of any $\xi \in \mathcal{S}(X_A)$. Assuming no preprocessing other than extracting and storing all the edges in the scene (and no other shortcuts in evaluating a test), this is roughly proportional to the actual algorithmic cost in CPU terms.

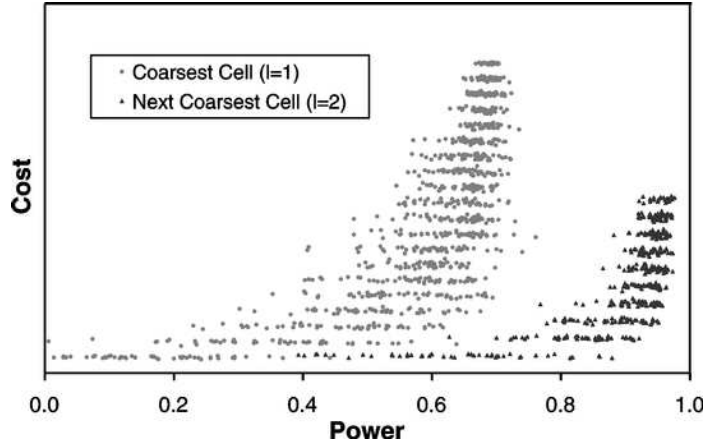

Fig. 8. *Cost vs. power curve for attributes of depth one and two.*

Recall our basic constraint: $P(X_A = 1 | Y \in A) = 1$ for every test $X_A$. In particular, we demand that $X_A = 1$ when the image data surrounding $R$ contains a rectangle whose pose belongs to $A$. Of course the test may also respond positively in the absence of such a rectangle, due to clutter and nearby rectangles; the likelihood of this happening is precisely $1 - \beta(A) = P(X_A = 1 | Y = 0)$. Intuitively, we expect that high power will only be possible at low invariance (specific poses). The power $\beta(A)$ is estimated from large samples of randomly selected background subimages.

In Figure 8 we plot cost versus power for the family of all tests generated for the root cell, $A_1$, referred to as "cell 1," and one of its two daughter cells, referred to as "cell 2." Thus each point is a pair $(\beta, c(X_{A,\beta}))$. For the root cell we cannot make tests with arbitrarily large power, at least not with such simple functionals. The "best tests" are those which are not strictly dominated by another test with respect to both cost and power—basically the convex envelope of the whole family; plots are given in [6]. Plots for cells at other depths are very similar, and the convexity assumption made in Sections 5 and 6 seems to be roughly satisfied.

Finally, one can ask whether the functional form of our global cost model, namely $c(X_{A,\beta}) = \Gamma(|A|) \times \Psi(\beta)$, is consistent with the data. This means an additive model for the log of the cost. In Figure 9 we plot the (base 2) logarithm of cost against the (base 2) logarithm of $\eta_l$ for five selected powers. Each point is one test—the one with lowest cost among those with power very close to a selected value. *The fact that the curves are roughly translations of each other is consistent with the additive model for the log-cost.* The roughly linear dependence of the log-cost with respect to $\log \Gamma(|A|)$ suggests a power dependence as a first approximation ($\Psi(x) \propto x^\alpha$ for some $\alpha \in [0, 1]$).

8.2.3. *Detection results.* We use the framework of Section 5—power-based cost for a fixed hierarchy. More specifically, from all the "best tests" created, we extracted one for each cell $A \in \mathcal{A}$ such that all the powers and costs are (approximately) the same at each level, which yields one sequence $(\beta_l, c_l)$, $l = 1, \ldots, 6$, which is increasing in both components and plotted in Figure 10 (left). Since the powers are increasing, the conditions of Corollary 2 are satisfied under the cost model. However, we need not assume that the cost model is valid; we can directly check whether $(\beta_l, c_l)$ satisfies the hypotheses of Corollary 1. In Figure 10 (right) we show, level by level, the (logarithms of the) values representing the two sides of (16). Clearly the conditions of Corollary 1 are easily satisfied.

The detection results for one scene are shown in Figure 11. In order to estimate total computation, we processed an $858 \times 626$ scene 100 times. The average time is 3.25 s on a Pentium 1.5 GHz. For comparison, we can perform an ideal hypothesis test for each fine cell ($Y \in A_{6,k}, k = 1, \ldots, 64$) based on simply counting *all* the edges in the region generated by the union of silhouettes over the poses in $A_{6,k}$ (a form of template-matching) and setting a threshold to obtain no false negatives. (This is a more discriminating test than $X_A$ for a fine cell $A$ because the latter uses only *some* of the edges.) The average processing time for this brute force approach is far larger (2338 s) but the results are virtually perfect. Finally, we can perform a two-stage analysis, first executing the CTF search and then doing the template-matching only at the detected poses. The processing time is virtually the same as for the

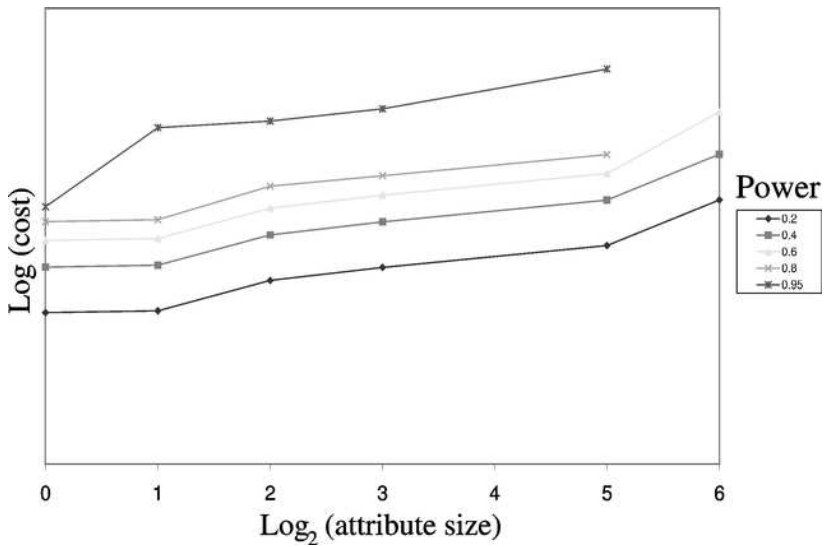

Fig. 9. *Log-cost vs. log-invariance for various powers.*

CTF search (about 3 s) but most of the false positives are removed; see Figure 11.

**9. Discussion and conclusion.** There are many problems in machine learning and perception which come down to differentiating among an enormous number of competing explanations, some very similar to each other and far too many to examine one-by-one. In these cases, efficient representations may be as important as statistical learning [18], and thinking about computation at the start of the day may be essential. It then seems prudent to model the computational process itself and hierarchical designs are a natural way to do this. Moreover, there is plenty of evidence that this works in practice. On the mathematical side, the questions that naturally arise from thinking about CTF representations and CTF search are of interest in themselves. We have provided one possible formulation; others could be envisioned.

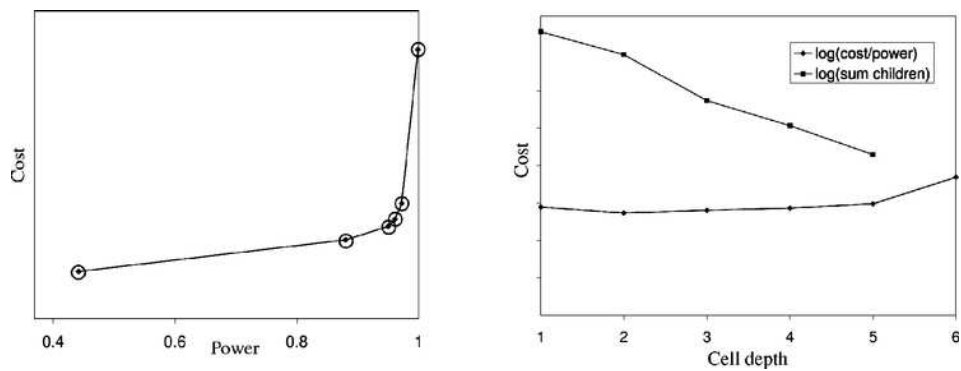

FIG. 10. Left: The pairs $(\beta_l, c_l)$ for the fixed hierarchy used in the experiments. Right: top curve: $l \to \log(C_l \times (c_{l+1}/\beta_{l+1}))$ where $C_l$ is the number of children of a node at level $l$; bottom curve: $l \to \log(c_l/\beta_l)$. The conditions of Corollary 1 are clearly satisfied.

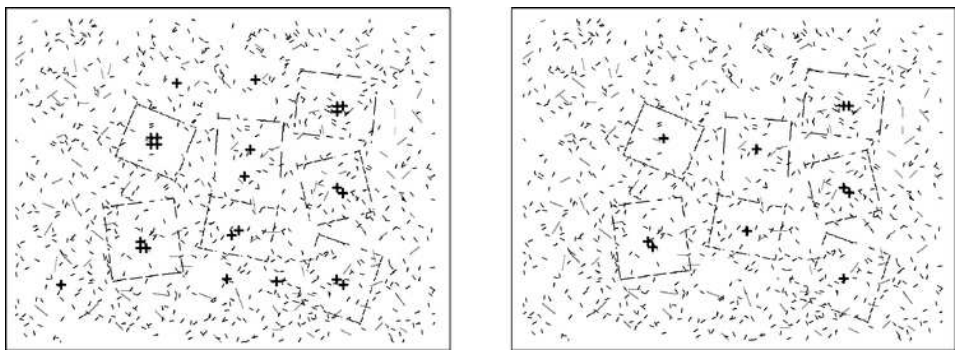

FIG. 11. Example of a detection result; the small crosses indicate the detected locations. Left: CTF detection only. Notice there are scattered false positives. Right: CTF search followed by template-matching. Nearly all the false positives are removed with virtually no increase in computation.

9.1. *Decision trees.* Of course "twenty questions," and the search strategies $T$ studied here based on a fixed family $\mathcal{X}$ of binary tests, invoke decision (or classification) trees—adaptive procedures for discriminating amongst hypotheses based on sequential testing.

Most of the literature is about an inductive framework. Trees are induced from a training set of i.i.d. samples from the joint distribution of the feature vector and class label, and binary tests result from comparing one component of the feature vector to a threshold. A tree is built in a top-down, greedy, recursive fashion based on some splitting criterion, usually entropy reduction [7]. The construction is then *data-driven* and *locally optimized*, guided by uncertainty reduction. There is a large literature on application of decision trees to pattern recognition which is outside the scope of this paper; see [1].

Generally, efficient (online) execution is not a criterion for construction or performance; for instance, the CART algorithm does not account for mean path length, let alone "costs" for the tests. Not surprisingly, recursive greedy designs are often globally inefficient, for instance in terms of the mean depth necessary to reach a given classification rate. A rarely studied alternative is to begin with an explicit statistical model for features and labels and compute a tree according to a global criterion involving *both* accuracy *and* (online) computation. The construction is then *model-driven* and *globally optimized*. Our approach to calculating $\widehat{Y}$ is of this general nature.

We refer the reader to [6] for an expanded discussion of these issues, including some early work due to Garey [16] on optimal testing procedures; related strategies for image retrieval [28]; comparisons between depth-first CTF and vanilla CART [20], showing that, in general, the latter is *not* CTF; and a special (if unrealistic) case, traced back to [10] and at the intersection of sequential statistics [8], game theory [5] and adaptive control processes [4], in which globally optimal testing strategies can be computed using dynamic programming, at least for "small problems." (See also work on cost-minimizing sequential procedures for Markov decision processes in [23].) In this special case, some comparisons in accuracy (resp. mean depth) between local and global strategies are given in [17] at a fixed mean depth (resp. accuracy), revealing an enormous difference in favor of global strategies, especially with skewed priors, that is, when a priori some classes are much more likely than others, which is precisely the situation in pattern recognition.

9.2. *Open issues.* Within our formulation there are some unanswered but fundamental mathematical questions and a few dubious assumptions. To begin with, we have divided the whole classification problem into two distinct and successive phases, first noncontextual (testing against nonspecific alternatives) and second contextual (testing one subset of explanations against another). We have shown that CTF search is effective, even optimal, in the first phase and preliminary results (not reported here) indicate the same is true of the second phase. However, whereas sensible, this division was artificially imposed; in particular, we have not shown that it emerges naturally from a global formulation of the problem. One might, for example, expand the family $\mathcal{X} = \{X_A, A \in \mathcal{A}\}$ into a much larger family of hypothesis tests for testing $Y \in A$ versus $Y \in B$ for various subsets $A, B$ and levels of error, and then attempt to *prove* that it is in fact computationally efficient to start with $B = A^c$ under some distributional assumptions, and reasonable trade-offs among scope, error and cost.

Whereas our results on fixed-powers hierarchies are fairly comprehensive, the results on variable-power hierarchies are evidently not. What is special, if anything, about the "harmonic cost function"? Simulations suggest that the CTF is generically optimal but we have not been able to prove this in general.

On the other hand, several of our model assumptions can be considered as too simplistic. Perhaps the cost model should be revisited; in simulations high power is not always attainable at high invariance (regardless of cost), at least for relatively simple tests (recall Figure 8). As pointed out earlier, supposing conditional independence under $P_0$ is disputable. Ideally, one should examine nontrivial dependency structures for $\mathcal{X}$, one appealing model being a first-order Markov structure of the tests as already depicted in the simulations of Section 5.4. Also, measuring computation under $P_0$ only is suspect. At some point in the computational process, as evidence accumulates from positive test results for the presence of a pattern of interest, the background hypothesis ceases to be dominant and all the class-conditional distributions must enter the story.

More ambitiously, an even more general optimization problem could be considered: Design the entire system including the subsets to be tested (not requiring a hierarchical structure a priori) as well as the levels of discrimination. This would likely involve a dependency structure for overlapping tests. Some of these questions are currently being investigated.

## APPENDIX

PROOF OF THEOREM 4. Consider a given tree-structured hierarchy $\mathcal{A}$. In this proof, we are mainly interested in the graph structure of $\mathcal{A}$. Here again it will be easier to consider the equivalent "augmented" model $\overline{\mathcal{A}}$ (see Section 4.5.4), thereby assuming the original $\mathcal{A}$ has been extended one level by adding a single child to each original leaf (in order to accommodate the perfect tests $X_{\{y\},1}$ which are performed at the end of the search for all $y \in \widehat{Y}$). Except for the power-one constraint for the final singleton tests: For any node $s$ in a strategy tree, the assigned power $\beta(s)$ may be freely chosen independently of how it is chosen when the corresponding attribute $A(s)$ appears at other nodes. Of course there must be no errors under $P_0$, but this is automatically satisfied by definition for any CTF strategy.

To prove the theorem we will proceed by recurrence over the size of subhierarchies of $\overline{\mathcal{A}}$. We actually need slightly more general objects than conventional subhierarchies (i.e., subtrees). We will call $\mathcal{H}$ a *generalized subhierarchy* if $\mathcal{H}$ is a finite union of subhierarchies of $\overline{\mathcal{A}}$. The cardinality of $\mathcal{H}$ is defined as the number of its nodes (internal or leaves). A CTF strategy for $\mathcal{H}$ satisfies the usual hypothesis that an attribute is tested if and only if all of its ancestors in $\mathcal{H}$ have been tested and returned a positive answer. Finally, for a node $B$ of $\overline{\mathcal{A}}$, denote by $\mathcal{H}_B$ the generalized subhierarchy composed of all strict descendents of $B$, in other words the union of all the subhierarchies rooted in direct children of $B$.

Now we prove by recurrence on the size $c$ of generalized subhierarchies which have the following property:

(P($c$))  *For any generalized subhierarchy $\mathcal{H}$ of $\overline{\mathcal{A}}$ of cardinality at most $c$, every CTF strategy with optimal choice of powers has the same cost $\mathcal{C}_{\mathrm{ctf}}(\mathcal{H})$. Furthermore, for any node $B \in \mathcal{H}$, the test $X_B$ is always performed in such a CTF strategy with the same power $\beta_B$, and this value depends only on $\mathcal{H}_B$, being therefore independent of the CTF strategy considered. Finally, if $\mathcal{H}$ is the union of several disjoint subhierarchies of $\overline{\mathcal{A}}$, then the CTF cost of $\mathcal{H}$ is the sum of the CTF costs of these subhierarchies.*

For $c = 1$, any generalized subhierarchy $\mathcal{H}$ must be a single node (attribute) corresponding to a perfect test, in which case the property is trivial.

Suppose (P($c$)) is true and consider a generalized subhierarchy $\mathcal{H}$ of cardinality $c+1$. Let $T$ be a CTF strategy for $\mathcal{H}$ with optimally chosen powers and let $B$ be the attribute which is tested at the root of $T$; necessarily $B$ has no ancestors in $\mathcal{H}$. Write $\overline{\mathcal{H}}_B$ for the generalized subhierarchy $\mathcal{H} \setminus (\{B\} \cup \mathcal{H}_B)$.

If $B$ is a leaf, then, by construction, its power is fixed to 1 and $\mathcal{H}_B = \varnothing$. Hence, after $B$ is tested with power 1 (thus returning a null answer under $P_0$), the remaining part of $T$ is a CTF strategy for subhierarchy $\overline{\mathcal{H}}_B$, and therefore, by the hypothesis of recurrence,

$$(28) \qquad E_0[C(T)] = \Psi(1) + \mathcal{C}_{\mathrm{ctf}}(\overline{\mathcal{H}}_B).$$

Suppose now that $B$ is not a leaf. If the test $X_B = 0$, the subsequent part of strategy $T$ must be a CTF exploration, with optimal powers, of the subhierarchy $\overline{\mathcal{H}}_B$. Similarly, if $X_B = 1$, the subsequent part of $T$ is a CTF strategy for $\overline{\mathcal{H}}_B \cup \mathcal{H}_B$, a disjoint union. By Lemma 3 and the recurrence hypothesis concerning cost additivity over disjoint subhierarchies, we therefore have

$$
\begin{aligned}
(29) \quad E_0[C(T)] &= \mathcal{C}_{\mathrm{ctf}}(\overline{\mathcal{H}}_B) + \Phi_{\Gamma(|B|)}(\mathcal{C}_{\mathrm{ctf}}(\overline{\mathcal{H}}_B \dot{\cup} \mathcal{H}_B) - \mathcal{C}_{\mathrm{ctf}}(\overline{\mathcal{H}}_B)) \\
&= \mathcal{C}_{\mathrm{ctf}}(\overline{\mathcal{H}}_B) + \Phi_{\Gamma(|B|)}(\mathcal{C}_{\mathrm{ctf}}(\mathcal{H}_B)).
\end{aligned}
$$

Furthermore, the second part of Lemma 3 shows that the optimal power chosen for $X_B$ only depends on $\mathcal{C}_{\mathrm{ctf}}(\mathcal{H}_B)$.

Property (P($c+1$)) is now an immediate consequence of (28) and (29), which concludes the proof. $\square$

PROOF OF THEOREM 6. Our goal is to prove the switching property (25) for the harmonic cost, that is,

$$
\forall\, y \geq x \geq 0, \forall\, a \geq b \geq 0 \qquad \Phi_a(x) + \Phi_b(\Phi_a(y) - \Phi_a(x)) \geq \Phi_a(x + \Phi_b(y - x)).
$$

This is obviously satisfied when $x = y$ (for any choice of $a$ and $b$). Denote by $C_L(y; x)$ [resp. $C_R(y; x)$] the left-hand (resp. right-hand) side of the above inequality. We will show that

$$
(30) \qquad \frac{\partial C_L(y; x)}{\partial y} \geq \frac{\partial C_R(y; x)}{\partial y} \qquad \text{for all } y \geq x \geq 0,
$$

which will conclude the proof. Taking derivatives in (30) we obtain

$$
(31) \qquad \Phi_b'(\Phi_a(y) - \Phi_a(x))\Phi_a'(y) \geq \Phi_a'(x + \Phi_b(y - x))\Phi_b'(y - x)
$$

with

$$
\Phi_b'(x) = \left( \frac{b}{x+b} \right)^2.
$$

After some elementary algebra, we find that (31) is equivalent to

$$
(y - x)[(x + a)^2 - a^2] \geq 0,
$$

which is true since $y \geq x$. This concludes the proof. $\square$

**Acknowledgment.** We are grateful to Franck Jung for performing the experiments on rectangle detection and we refer the reader to his cited work for a convincing example of the "real thing"—the difficult task of detecting buildings in aerial photographs.

## REFERENCES

[1] AMIT, Y. (2002). 2*D Object Detection and Recognition*. MIT Press, Cambridge, MA.

[2] AMIT, Y. and GEMAN, D. (1999). A computational model for visual selection. *Neural Computation* **11** 1691–1715.

[3] AMIT, Y., GEMAN, D. and FAN, X. (2004). A coarse-to-fine strategy for multiclass shape detection. *IEEE Trans. Pattern Analysis and Machine Intelligence* **26** 1606–1621. MR2005782

[4] BELLMAN, R. (1961). *Adaptive Control Processes*: *A Guided Tour*. Princeton Univ. Press. MR134403

[5] BLACKWELL, D. and GIRSCHICK, M. A. (1954). *Theory of Games and Statistical Decisions*. Wiley, New York. MR70134

[6] BLANCHARD, G. and GEMAN, D. (2003). Hierarchical testing designs for pattern recognition. Technical Report 2003-07, Département de Mathématiques, Univ. de Paris-Sud.

[7] BREIMAN, L., FRIEDMAN, J., OLSHEN, R. and STONE, C. (1984). *Classification and Regression Trees*. Wadsworth, Belmont, CA. MR726392

[8] CHERNOFF, H. (1972). *Sequential Analysis and Optimal Design*. SIAM, Philadelphia. MR309249

[9] COOTES, T. F. and TAYLOR, C. J. (1996). Locating faces using statistical feature detectors. In *Proc. Second International Conference on Automatic Face and Gesture Recognition* 204–209. IEEE Press, New York.

[10] DEGROOT, M. H. (1970). *Optimal Statistical Decisions*. McGraw–Hill, New York. MR356303

[11] DESIMONE, R., MILLER, E. K., CHELAZZI, L. and LUESCHOW, A. (1995). Multiple memory systems in visual cortex. In *The Cognitive Neurosciences* (M. S. Gazzaniga, ed.) 475–486. MIT Press, Cambridge, MA.

[12] DIETTERICH, T. (2000). The divide-and-conquer manifesto. In *Proc. Eleventh International Conference on Algorithmic Learning Theory. Lecture Notes in Artificial Intelligence* **1968** 13–26. Springer, New York.

[13] FLEURET, F. (2000). Détection hiérarchique de visages par apprentissage statistique. Ph.D. dissertation, Univ. Paris VI, Jussieu.

[14] FLEURET, F. and GEMAN, D. (2001). Coarse-to-fine face detection. *International J. Computer Vision* **41** 85–107.

[15] FRITSCH, J. and FINKE, M. (1998). Applying divide and conquer to large scale pattern recognition tasks. *Neural Networks*: *Tricks of the Trade. Lecture Notes in Comput. Sci.* (G. B. Orr and K.-R. Müller, eds.) **1524** 315–342. Springer, New York.

[16] GAREY, M. R. (1972). Optimal binary identification procedures. *SIAM J. Appl. Math.* **23** 173–186. MR401334

[17] GEMAN, D. and JEDYNAK, B. (2001). Model-based classification trees. *IEEE Trans. Inform. Theory* **47** 1075–1082. MR1829333

[18] GEMAN, S., BIENENSTOCK, E. and DOURSAT, R. (1992). Neural networks and the bias/variance dilemma. *Neural Computation* **4** 1–58.

[19] GEMAN, S., MANBECK, K. and MCCLURE, D. (1995). Coarse-to-fine search and rank-sum statistics in object recognition. Technical report, Div. Applied Mathematics, Brown Univ.

[20] JUNG, F. (2001). Reconnaissance d'objets par focalisation et detection de changements. Ph.D. dissertation, Ecole Polytechnique, Paris.

[21] JUNG, F. (2002). Detecting new buildings from time-varying aerial stereo pairs. Technical report, IGN.

[22] Osuna, E., Freund, R. and Girosit, F. (1997). Training support vector machines: An application to face detection. In *Proc. Computer Vision and Pattern Recognition* 130–136. IEEE Press, New York.

[23] Puterman, M. (1994). *Markov Decision Processes.* Wiley, New York. MR1270015

[24] Rowley, H. A., Baluja, S. and Kanade, T. (1998). Neural network-based face detection. *IEEE Trans. Pattern Analysis and Machine Intelligence* **20** 23–38.

[25] Socolinsky, D. A., Neuheisel, J. D., Priebe, C. E., De Vinney, J. and Marchette, D. (2002). Fast face detection with a boosted cccd classifier. Technical report, Johns Hopkins Univ.

[26] Sung, K.-K. and Poggio, T. (1998). Example-based learning for view-based human face detection. *IEEE Trans. Pattern Analysis and Machine Intelligence* **20** 39–51.

[27] Thorpe, S., Fize, D. and Marlot, C. (1996). Speed of processing in the human visual system. *Nature* **381** 520–522.

[28] Trouvé, A. and Yu, Y. (2002). Entropy reduction strategies on tree structured retrieval spaces. In *Proc. Colloquium on Mathematics and Computer Science. II. Algorithms, Trees, Combinatorics and Probabilities* (B. Chauvin, P. Flajolet, D. Gardy and A. Mokkadem, eds.) 513–525. Birkhäuser, Basel. MR1940157

[29] Viola, P. and Jones, M. J. (2001). Robust real-time face detection. In *Proc. Eighth IEEE International Conference on Computer Vision* **2** 747. IEEE Press, New York.

[30] Yuille, A. L., Hallinan, P. and Cohen, D. S. (1992). Feature extraction from faces using deformable templates. *International J. Computer Vision* **8** 99–111.

CNRS
France
and
IDA Group
Fraunhofer FIRST
Berlin
Germany
E-mail: Gilles.Blanchard@first.fraunhofer.de

Department of Applied Mathematics
  and Statistics
Johns Hopkins University
Baltimore, Maryland 21218-2682
USA
E-mail: geman@jhu.edu